\documentclass[11pt]{article} 
 
\usepackage{latexsym} 
\usepackage{amsfonts} 
\usepackage{amsmath} 
\usepackage{amsthm} 
\usepackage{mathrsfs} 
\usepackage{dsfont} 
\usepackage{bbold} 
\usepackage[english]{babel} 
\usepackage{caption} 
\usepackage{epsfig} 
\usepackage{float} 
\usepackage{subfigure} 
\usepackage{psfrag} 
\usepackage{graphicx} 
\usepackage{epsfig}

\DeclareMathOperator{\Val}{\matV}

\newtheorem{theorem}{Theorem} 
\newtheorem*{prop*}{Theorem}

\newtheorem{defi}[theorem]{Definition} 
\newtheorem{lemma}[theorem]{Lemma} 
 
\newtheorem{rmk}[theorem]{Remark}

\newtheorem{hyp}{Hypothesis} 
 
\newcommand{\zerarcounters}{\setcounter{equation}{0}\setcounter{theorem}{0}} 
 
\newcommand{\ZZZ}{\mathds{Z}} 
 
\newcommand{\NNN}{\mathds{N}} 
 
\newcommand{\RRR}{\mathds{R}} 
\newcommand{\TTT}{\mathds{T}}

\newcommand{\calA}{{\mathcal A}} 
\newcommand{\BB}{{\mathcal B}} 
\newcommand{\CCCC}{{\mathcal C}}

\newcommand{\calF}{{\mathcal F}} 
\newcommand{\calG}{{\mathcal G}}

\newcommand{\MM}{{\mathcal M}}

\newcommand{\calP}{{\mathcal P}} 
 
\newcommand{\RR}{{\mathcal R}}

\newcommand{\gotn}{{\mathfrak n}}

\newcommand{\gotF}{{\mathfrak F}}

\newcommand{\gotL}{{\mathfrak L}} 
 
\newcommand{\gotN}{{\mathfrak N}}

\newcommand{\gotR}{{\mathfrak R}} 
\newcommand{\gotS}{{\mathfrak S}}

\newcommand{\matV}{{\mathscr V}}

\newcommand{\ol}{\overline} 
 
\newcommand{\Fullbox}{{\rule{2.0mm}{2.0mm}}} 
 
\newcommand{\EP}{\hfill\Fullbox\vspace{0.2cm}} 
\newcommand{\prova}{\noindent{\it Proof. }} 
\newcommand{\io}{\infty} 
\newcommand{\e}{\varepsilon} 
\newcommand{\al}{\alpha} 
 
\newcommand{\be}{\beta} 
\newcommand{\n}{\nu} 
 
\newcommand{\x}{\xi} 
\newcommand{\p}{\pi} 
 
\newcommand{\g}{\gamma}

\newcommand{\la}{\lambda} 
 
\newcommand{\s}{\sigma}

\newcommand{\oo}{\boldsymbol{\omega}} 
 
\newcommand{\aaa}{\boldsymbol{\alpha}}

\newcommand{\nn}{\boldsymbol{\nu}} 
\newcommand{\pps}{\boldsymbol{\psi}} 
\newcommand{\vzero}{\boldsymbol{0}}

\newcommand{\AAA}{\boldsymbol{A}} 
 
\newcommand{\ii}{{\rm i}}

\def\ins#1#2#3{\vbox to0pt{\kern-#2 \hbox{\kern#1 #3}\vss}\nointerlineskip} 
 
\begin{document} 
 
\title{\bf Response solutions for arbitrary quasi-periodic 
perturbations with Bryuno frequency vector} 
 
\author 
{\bf Livia Corsi and Guido Gentile 
\vspace{2mm} 
\\ \small 
Dipartimento di Matematica, Universit\`a di Roma Tre, Roma, 
I-00146, Italy 
\\ \small 
E-mail: lcorsi@mat.uniroma3.it, gentile@mat.uniroma3.it} 
 
\date{} 
 
\maketitle 
 
\begin{abstract} 
We study the problem of existence of response solutions for a 
real-analytic one-dimensional system, consisting of a rotator 
subjected to a small quasi-periodic forcing. We prove that at least 
one response solution always exists, without any assumption 
on the forcing besides smallness and analyticity. 
This strengthens the results available in the literature, 
where generic non-degeneracy conditions are assumed. 
The proof is based on a diagrammatic formalism and relies on 
renormalisation group techniques, which exploit the formal analogy 
with problems of quantum field theory; a crucial role is played by 
remarkable identities between classes of diagrams.
\end{abstract} 
 
 
 
 
\zerarcounters 
\section{Introduction} 
\label{sec:1} 
 
Consider the one-dimensional system 
\begin{equation}\label{eq:1.1} 
\ddot \be = -\e F(\oo t, \be),\qquad 
F(\oo t, \be):= \partial_{\be}f(\oo t,\be) , 
\end{equation} 
where $\be\in\TTT=\RRR/2\p\ZZZ$, $f\!:\TTT^{d+1}\to\RRR$ is a real-analytic 
function, $\oo\in\RRR^{d}$ and $\e$ is a real number, called the 
\emph{perturbation parameter}; 
hence the \emph{forcing function} (or perturbation) $F$ is 
quasi-periodic in $t$, with \emph{frequency vector} $\oo$. 
 
It is well known that, for $d=1$ (periodic forcing) and $\e$ small enough, 
there exist periodic solutions to (\ref{eq:1.1}) with the same period as the 
forcing. In fact the existence of periodic solutions to (\ref{eq:1.1}), 
or to the more general equation 
\begin{equation} \label{eq:1.2} 
\ddot \be = -\partial_{\be}V(\be)-\e F(\oo t, \be), 
\end{equation} 
with $V\!:\RRR\to\RRR$ real-analytic, can be discussed by relying on 
\emph{Melnikov method} \cite{CH,GH}. 
A possible approach consists in splitting the equations of motion into 
two separate equations, the so-called \emph{range equation} and 
\emph{bifurcation equation}. 
Then, one can solve the first equation in terms of a free parameter, 
and then fix the latter by solving the second equation (which 
represents an implicit function problem). This is usually done by 
assuming some \emph{non-degeneracy condition} involving the
perturbation, and this entails the analyticity  of the solution.
If no such condition is assumed, a result of the same kind still 
holds \cite{ZL,ALGM,CG}, but the scenario appears slightly more 
complicated: for instance the persisting periodic solutions are no 
longer analytic in the perturbation parameter. 
 
If the forcing is quasi-periodic, one can still study the problem of 
existence of quasi-periodic solutions with the same frequency vector $\oo$ 
as the forcing, for $\e$ small enough. The analysis becomes much more 
involved, because of the small divisor problem. However, under 
some generic non-degeneracy condition, the analysis can be carried 
out in a similar way and the bifurcation scenario can be described 
in a rather detailed way; see for instance \cite{BHY}. On the contrary, 
if no assumption at all is made on the perturbation, the small divisor 
problem and the implicit function problem become inevitably tangled 
together and new difficulties arise. In this paper we focus on this 
situation, so we study (\ref{eq:1.1}) without making any assumption 
on the forcing function besides analyticity. Of course, we shall make some
assumption 
of strong irrationality on the frequency vector $\oo$, say we shall 
assume some mild Diophantine condition, such as the Bryuno condition 
(see below). 
 
Note that (\ref{eq:1.1}) can be seen as the Hamilton equations for the 
system described by the Hamiltonian function 
\begin{equation}\label{eq:1.3} 
H(\aaa,\be,\AAA,B)=\oo\cdot \AAA + \frac{1}{2}B^{2}+\e f(\aaa,\be), 
\end{equation} 
where $\oo \in \RRR^{d}$ is fixed, $(\aaa,\be)\in \TTT^{d}\times \TTT$ 
and $(\AAA,B)\in \RRR^{d}\times \RRR$ are conjugate variables and 
$f$ is an analytic periodic function of $(\aaa,\be)$. Indeed, the 
corresponding Hamilton equations for the angle variables are closed, 
and are given by 
\begin{equation}\label{eq:1.4} 
\dot \aaa = \oo,\qquad 
\ddot \be = -\e \partial_{\be}f(\aaa,\be), 
\end{equation} 
that we can rewrite as (\ref{eq:1.1}). 
Therefore the problem of existence of \emph{response solutions}, i.e. 
quasi-periodic solutions to (\ref{eq:1.1}) with frequency vector $\oo$, 
can be seen as a problem of persistence of lower-dimensional 
(or resonant) tori, more precisely of $d$-dimensional tori 
for a system with $d+1$ degrees of freedom. In the case (\ref{eq:1.3}) 
the \emph{unperturbed} (i.e. with $\e=0$) Hamiltonian is isochronous in all
but one 
angle variables. The existence of $d$-dimensional tori in systems 
with $d+1$ degrees of freedom, without imposing any 
non-degeneracy condition on the perturbation except analyticity, 
was first studied by Cheng \cite{Ch}. He proved that, 
for convex unperturbed Hamiltonians, there exists at least one 
$d$-dimensional torus continuing a $d$-dimensional submanifold 
of the $d+1$ unperturbed resonant torus on which the flow 
is quasi-periodic with frequency vector $\oo\in\RRR^{d}$ 
satisfying the standard Diophantine condition 
$|\oo\cdot\nn|\ge \g|\nn|^{-\tau}$ for all 
$\nn\in\ZZZ^{d}\setminus\{\vzero\}$, 
and for some $\g>0$ and $\tau>d-1$ (here and henceforth 
$\cdot$ denotes the standard scalar product in $\RRR^{d}$ 
and $|\nn|=|\nn|_{1}=|\n_{1}|+\ldots+|\n_{d}|$). 
 
We prove a result of the same kind for the equation (\ref{eq:1.1}), 
that is the existence of at least one response solution for $\e$ 
small enough -- see Theorem 
\ref{thm:2.2} in Section \ref{sec:2}. 
Even if the system (\ref{eq:1.3}) can be seen as a simplified model for the 
problem of lower-dimensional tori, we think that our result can be of interest 
by its own. First of all, Cheng's result does not directly apply, 
since both the convexity property he requires is obviously not satisfied 
by the Hamiltonian (\ref{eq:1.3}) and we allow a weaker Diophantine 
condition on the frequency vector. Moreover, 
just because of its simplicity, the model is particularly suited 
to point out the main issues of the proof, avoiding all 
aspects that would add only technical intricacies without 
shedding further light on the problem. 
Finally, our method is completely 
different: it is based on the analysis and resummation of the 
perturbation series through renormalisation group techniques, 
and not on an iteration scheme \emph{\`a la} KAM. In particular 
a crucial role in the proof will be played by remarkable identities 
between classes of diagrams. By exploiting the analogy of the method
with the techniques of quantum field theory, one can see
the solution as the one-point Schwinger function of a suitable
Euclidean field theory -- this has been explicitly shown in the
case of KAM tori \cite{GGM} --; then the identities between diagrams
can be imagined as due to a suitable Ward identity
that follows from the symmetries of the field theory --
again, this has been checked for the KAM theorem \cite{BGK},
and we leave it as a conjecture in our case.
 
\zerarcounters 
\section{Results} 
\label{sec:2} 
 
Consider equation (\ref{eq:1.1}) and take the solution for the 
unperturbed system given by $\be(t)=\be_{0}$. We want to study whether 
for some value of $\be_{0}$ such a solution can be continued under
perturbation. 
 
\begin{hyp}\label{hyp1} 
$\oo$ satisfies the Bryuno condition $\BB(\oo)<\io$, where 
\begin{equation} \nonumber 
\BB(\oo):=\sum_{m=0}^{\io}\frac{1}{2^{m}}\log\frac{1}{\al_{m}(\oo)}, 
\qquad 
\al_{m}(\oo):=\inf_{0<|\nn|\le 2^{m}}|\oo\cdot\nn|. 
\end{equation} 
\end{hyp} 
 
Write 
\begin{equation}\label{eq:2.1} 
f(\aaa,\be)=\sum_{\nn\in\ZZZ^{d}}f_{\nn}(\be)e^{\ii\nn\cdot\aaa}, \qquad 
F(\aaa,\be)=\sum_{\nn\in\ZZZ^{d}}F_{\nn}(\be)e^{\ii\nn\cdot\aaa}. 
\end{equation} 
%
 
\begin{hyp}\label{hyp2} 
$\be_{0}^{*}$ is a zero of order $\gotn$ for $F_{\vzero}(\be)$ with 
$\gotn$ odd. Assume also $\e \partial_{\be}^{\gotn}F_{\vzero}(\be_{0}^{*})<0$ 
for fixed $\e\neq0$. 
\end{hyp} 
 
Eventually we shall want to get rid of 
Hypothesis \ref{hyp2}: however, we shall first assume it 
to simplify the analysis, and at the end we shall show how to remove it. 
 
We look for a solution to (\ref{eq:1.1}) of the form $\be(t)=\be_{0}+b(t)$, 
with 
\begin{equation} \label{eq:2.2} 
b(t)=\sum_{\nn\in\ZZZ^{d}_{*}}e^{\ii\nn\cdot\oo t}b_{\nn} 
\end{equation} 
where $\ZZZ^{d}_{*}=\ZZZ^{d}\setminus\{\vzero\}$. 
In Fourier space (\ref{eq:1.1}) becomes 
\begin{subequations} 
\begin{align} 
& (\oo\cdot\nn)^{2}b_{\nn}=\e[F(\oo t,\be)]_{\nn}, \quad \nn\neq\vzero,
\label{eq:2.3a} \\ 
& [F(\oo t,\be)]_{\vzero}=0 , 
\label{eq:2.3b} 
\end{align} 
\label{eq:2.3} 
\end{subequations} 
where 
\begin{equation} \nonumber 
[F(\pps,\be)]_{\nn} = 
\sum_{r\ge 0} \sum_{\substack{\nn_{0}+\ldots+\nn_{r}=\nn \\ 
\nn_{0}\in \ZZZ^{d} \\ 
\nn_{i} \in \ZZZ^{d}_{*},\, i=1,\ldots,r}} 
\frac{1}{r!}\partial^{r}_{\be} F_{\nn_{0}}(\be_{0})\,\prod_{i=1}^{r}b_{\nn_{i}} . 
\end{equation}  

Our first result will be the following. 
 
\begin{theorem}\label{thm:2.1} 
Consider the equation (\ref{eq:1.1}) and assume Hypotheses \ref{hyp1} and 
\ref{hyp2}. If $\e$ is small enough, there exists at least one 
quasi-periodic solution $\be(t)$ to (\ref{eq:1.1}) with frequency vector 
$\oo$,  such that $\be(t) \to\be_{0}^{*}$ as $\e\to0$. 
\end{theorem} 
 
The proof will be carried out through Sections \ref{sec:3} to \ref{sec:5}. 
First, after introducing the basic notations in Section \ref{sec:3}, 
we shall show in Section \ref{sec:4} that, 
under the assumption that further conditions are satisfied, 
for $\e$ small enough and arbitrary $\be_0$ there exists a solution 
\begin{equation}\label{eq:2.4} 
\be(t)=\be_{0}+b(t;\e,\be_{0}), 
\end{equation} 
to (\ref{eq:2.3a}), depending on $\e,\be_{0}$, with $b(t)=b(t;\e,\be_{0})$ 
a zero-average function. For such a solution define 
\begin{equation}\label{eq:2.5} 
G(\e,\be_{0}):=[F(\oo t,\be(t))]_{\vzero}, 
\end{equation} 
and consider the implicit function equation 
\begin{equation}\label{eq:2.6} 
G(\e,\be_{0})=0. 
\end{equation} 
Then we shall prove in Section \ref{sec:5} that one can fix $\be_{0}= 
\be_{0}(\e)$ in such a way that (\ref{eq:2.6}) holds and the conditions 
mentioned above are also satisfied. Hence for that 
$\be_{0}(\e)$ the function (\ref{eq:2.4}) is a solution of the whole 
system (\ref{eq:2.3}). 
 
Next, we shall see how to remove Hypothesis \ref{hyp2} in order to prove 
the existence of a response solution without any assumption on the forcing 
function, so as to obtain the following result, which is the main result of 
the paper. 
 
\begin{theorem}\label{thm:2.2} 
Consider the equation (\ref{eq:1.1}) and assume Hypothesis \ref{hyp1}. 
There exists $\e_{0}>0$ such that for all $\e$ with $|\e|<\e_{0}$ there 
is at least one quasi-periodic solution to (\ref{eq:1.1}) with frequency 
vector $\oo$. 
\end{theorem} 
 
Note that if $F_{\vzero}(\be)$ does not identically vanish, then 
Theorem \ref{thm:2.2} follows immediately from Theorem \ref{thm:2.1}. 
Indeed, the function $f_{\vzero}(\be)$ is analytic and periodic, 
hence, if it is not identically constant, it has at least one maximum 
point $\be_{0}'$ and one minimum point $\be_{0}''$, where 
$\partial_{\be}^{\gotn'+1}f_{\vzero}(\be_{0}')<0$ 
and $\partial_{\be}^{\gotn''+1}f_{\vzero}(\be_{0}'')>0$, 
for some $\gotn'$ and $\gotn''$ both odd. 
Let $\e$ be fixed small enough, say $|\e|<\e_{0}$ for a suitable 
$\e_{0}$: choose $\be_{0}^{*}=\be_{0}'$ if $\e>0$ and 
$\be_{0}^{*}=\be_{0}''$ if $\e<0$. Then Hypothesis \ref{hyp2} is 
satisfied, and we can apply Theorem \ref{thm:2.1} to deduce the 
existence of a quasi-periodic solution with frequency vector $\oo$. 
However, the function $f_{\vzero}(\be)$ can be identically constant, 
and hence $F_{\vzero}(\be)$ can vanish identically, so that some further 
work will be needed to prove Theorem \ref{thm:2.2}: this will be 
performed in Section \ref{sec:6}. 
 
\zerarcounters 
\section{Diagrammatic rules and multiscale analysis} 
\label{sec:3} 
 
We want to study whether it is possible to express the function 
$b(t;\e,\be_{0})$ appearing in (\ref{eq:2.4}) 
as a convergent series. Let us start by writing formally 
\begin{equation}\label{eq:3.1} 
b(t;\e,\be_{0})=\sum_{k\ge 1}\e^{k}b^{(k)}(t;\be_{0})= 
\sum_{k\ge 1}\e^{k}\sum_{\nn\in\ZZZ^{d}_{*}} 
e^{\ii\nn\cdot\oo t}b^{(k)}_{\nn}(\be_{0}). 
\end{equation} 
If we define recursively for $k\ge 1$ 
\begin{equation}\label{eq:3.2} 
b_{\nn}^{(k)}(\be_{0})=\frac{1}{(\oo\cdot\nn)^{2}} 
[F(\oo t,\be)]^{(k-1)}_{\nn}, 
\end{equation} 
where $[F(\oo t,\be)]_{\nn}^{(0)}=F_{\nn}(\be_{0})$ and, for $k\ge1$, 
\begin{equation}\label{eq:3.3} 
[F(\oo t,\be)]^{(k)}_{\nn}=\sum_{s\ge 1} 
\sum_{\substack{\nn_{0}+\ldots+\nn_{s}=\nn \\ \nn_{0}\in \ZZZ^{d} \\ 
\nn_{i} \in \ZZZ^{d}_{*},\, i=1,\ldots,s}} 
\frac{1}{s!}\partial^{s}_{\be} 
F_{\nn_{0}}(\be_{0}) 
\sum_{\substack{k_{1}+\ldots+k_{s}=k, \\ k_{i}\ge 1}} 
\prod_{i=1}^{s} b_{\nn_{i}}^{(k_{i})}(\be_{0}) , 
\end{equation} 
the series (\ref{eq:3.1}) turns out to be a formal solution of (\ref{eq:2.3a}): 
the coefficients $b^{(k)}_{\nn}(\be_{0})$ are well defined for all 
$k\ge1$ and all $\nn\in\ZZZ^{d}_{*}$ -- by Hypothesis \ref{hyp1} -- 
and solve (\ref{eq:2.3a}) order by order -- as it is straightforward to check. 
 
Write also, again formally, 
\begin{equation}\label{eq:3.4} 
G(\e,\be_{0})=\sum_{k \ge 0}\e^{k}G^{(k)}(\be_{0}), 
\end{equation} 
with $G^{(0)}(\be_{0})=F_{\vzero}(\be_{0})$ and, for $k\ge1$ 
\begin{equation}\label{eq:3.5} 
G^{(k)}(\be_{0})=\sum_{s\ge 1} 
\sum_{\substack{\nn_{0}+\ldots+\nn_{s}=\vzero \\ \nn_{0}\in \ZZZ^{d}\\ 
\nn_{i} \in \ZZZ^{d}_{*}, \, i=1,\ldots,s}} 
\frac{1}{s!}\partial^{s}_{\be}F_{\nn_{0}}(\be_{0}) 
\sum_{\substack{k_{1}+\ldots+k_{s}=k, \\ k_{i}\ge 1}} 
\prod_{i=1}^{s} b_{\nn_{i}}^{(k_{i})}(\be_{0}) . 
\end{equation} 
Of course, Hypothesis \ref{hyp1} yields that the formal series (\ref{eq:3.4})
is  well-defined too. 
 
Unfortunately the power series (\ref{eq:3.1}) and (\ref{eq:3.4}) may not be 
convergent (as far as we know). However we shall see how to 
construct two series (convergent if $\be_{0}$ is suitably chosen)
whose formal expansion coincide with (\ref{eq:3.1}) and
(\ref{eq:3.4}). As we shall see, this leads to express the response
solution as a series of contributions each of which can be graphically
represented as a suitable diagram.

A graph is a set of points and lines connecting them. 
A \emph{tree} $\theta$ is a graph with no cycle, 
such that all the lines are oriented toward a unique point 
(\emph{root}) which has only one incident line $\ell_{\theta}$ 
(\emph{root line}). 
All the points in a tree except the root are called \emph{nodes}. 
The orientation of the lines in a tree induces a partial ordering 
relation ($\preceq$) between the nodes and the lines: we can 
imagine that each line carries an arrow pointing toward the root. 
Given two nodes $v$ and $w$, 
we shall write $w \prec v$ every time $v$ is along the path 
(of lines) which connects $w$ to the root. 
 
We denote by $N(\theta)$ and $L(\theta)$ the sets of nodes and 
lines in $\theta$ respectively. 
Since a line $\ell\in L(\theta)$ is uniquely identified 
with the node $v$ which it leaves, we may write $\ell = \ell_{v}$. 
We write $\ell_{w} \prec \ell_{v}$ if $w\prec v$, and $w\prec\ell=\ell_{v}$
if $w\preceq v$; 
if $\ell$ and $\ell'$ are two comparable lines, i.e. 
$\ell' \prec \ell$, we denote by $\calP(\ell,\ell')$ the 
(unique) path of lines connecting $\ell'$ to $\ell$, with $\ell$ and 
$\ell'$ not included (in particular $\calP(\ell,\ell')=\emptyset$ 
if $\ell'$ enters the node $\ell$ exits). 
 
With each node $v\in N(\theta)$ we associate a \emph{mode} label 
$\nn_{v}\in \ZZZ^{d}$ and we denote by $s_{v}$ the number of lines 
entering $v$. With each line $\ell$ 
we associate a \emph{momentum} $\nn_{\ell}\in \ZZZ^{d}_{*}$, 
except for the root line which can have either zero momentum or not, 
i.e. $\nn_{\ell_{\theta}}\in\ZZZ^{d}$. 
Finally, we associate with each line $\ell$ also a \emph{scale label} 
such that $n_{\ell}=-1$ if $\nn_{\ell}=\vzero$, while 
$n_{\ell}\in\ZZZ_{+}$ if $\nn_{\ell}\neq\vzero$. 
Note that one can have $n_{\ell}=-1$ only if $\ell$ is the root line 
of $\theta$. 
 
We force the following \emph{conservation law} 
\begin{equation}\label{eq:3.6} 
\nn_{\ell}=\sum_{\substack{ w \in N(\theta) \\ w\prec \ell}}\nn_{w}. 
\end{equation} 

In the following we shall call trees tout court the trees 
with labels, and we shall use the term \emph{unlabelled tree} 
for the trees without labels. 
 
We shall say that two trees are \emph{equivalent} if they can be 
transformed into each other by continuously deforming the lines in 
such a way that these do not cross each other and also labels match. 
This provides an equivalence relation on the set of the trees -- as it is 
easy to check. From now on we shall call trees tout court such 
equivalence classes. 
 
Given a tree $\theta$ we call \emph{order} of $\theta$ the 
number $k(\theta)=|N(\theta)|=|L(\theta)|$ (for any finite set $S$ 
we denote by $|S|$ its cardinality) and \emph{total momentum} of 
$\theta$ the momentum associated with $\ell_{\theta}$. 
We shall denote by $\Theta_{k,\nn}$ the set of  trees 
with order $k$ and total momentum $\nn$. More generally, 
if $T$ is a subgraph of $\theta$ (i.e. a set of nodes $N(T)\subseteq
N(\theta)$ connected by lines $L(T)\subseteq L(\theta)$), we call \emph{order}
of $T$ the number $k(T)=|N(T)|$. We say that a line enters $T$ if it 
connects a node $v\notin N(T)$ to a node $w\in N(T)$, 
and we say  that a line exits $T$ if it connects a node $v\in N(T)$ 
to a node $w\notin N(T)$. Of course, if a line $\ell$ enters or exits $T$,
then $\ell\notin L(T)$
 
\begin{rmk}\label{rmk:3.1a} 
\emph{ 
One has $\displaystyle{\sum_{v\in N(\theta)}s_{v}=k(\theta)-1}$. 
} 
\end{rmk} 
 
A \emph{cluster} $T$ on scale $n$ is a maximal subgraph 
of a tree $\theta$ such that all the lines have scales 
$n'\le n$ and there is at least a line with scale $n$. 
The lines entering the cluster $T$ and the line coming 
out from it (unique if existing at all) are called the 
\emph{external} lines of $T$. 
 
A \emph{self-energy cluster} is a cluster $T$ such that 
(i) $T$ has only one entering line $\ell'_{T}$ and 
one exiting line $\ell_{T}$, (ii) one has $\nn_{\ell_{T}}= 
\nn_{\ell'_{T}}$ and hence 
\begin{equation}\label{eq:3.7} 
\sum_{v\in N(T)}\nn_{v}=\vzero. 
\end{equation} 

For any self-energy cluster $T$, set $\calP_{T}= 
\calP(\ell_{T},\ell'_{T})$. More generally, if $T$ is a subgraph 
of $\theta$ with only one entering line $\ell'$ and one exiting line 
$\ell$, we can set $\calP_{T}=\calP(\ell,\ell')$. 
We shall say that a self-energy cluster is on 
scale $-1$, if $N(T)=\{v\}$ with of course $\nn_{v}=\vzero$ 
(so that $\calP_{T}=\emptyset$). 
 
A \emph{left-fake cluster} $T$ on scale $n$ is a connected 
subgraph of a tree $\theta$ with only one entering line $\ell'_{T}$ 
and one exiting line $\ell_{T}$ such that (i) all the lines in $T$ have 
scale $\le n$ and there is in $T$ at least a line on scale $n$, 
(ii) $\ell'_{T}$ is on scale $n+1$ 
and $\ell_{T}$ is on scale $n$, and (iii) one has 
$\nn_{\ell_{T}}=\nn_{\ell'_{T}}$. 
Analogously a \emph{right-fake cluster} $T$ on scale $n$ is 
a connected subgraph of a tree $\theta$ with only one entering 
line $\ell'_{T}$ and one exiting line $\ell_{T}$ such that (i) all the 
lines in $T$ have scale $\le n$ and there is in $T$ at least a line 
on scale $n$, (ii)  $\ell'_{T}$ is on scale $n$ 
and $\ell_{T}$ is on scale $n+1$, and (iii) one has 
$\nn_{\ell_{T}}=\nn_{\ell'_{T}}$. 
Roughly speaking, a left-fake (respectively right-fake) cluster $T$ fails to 
be a self-energy cluster only because the exiting (respectively the 
entering) line is on scale equal to the scale of $T$. 
 
\begin{rmk}\label{rmk:3.1} 
\emph{ 
Given a self-energy cluster $T$, the momenta of the lines in $\calP_{T}$ 
depend on $\nn_{\ell'_{T}}$ because of the conservation law (\ref{eq:3.6}). 
More precisely, for all $\ell\in\calP_{T}$ one has 
$\nn_{\ell}=\nn_{\ell}^{0}+\nn_{\ell'_{T}}$ with 
\begin{equation}\label{eq:3.8} 
\nn_{\ell}^{0}=\sum_{\substack{w\in N(T) \\ w\prec \ell}} 
\nn_{w},
\end{equation} 
while all the other labels in $T$ do not depend on $\nn_{\ell'_{T}}$. Clearly,
this holds also for left-fake and right-fake clusters.
} 
\end{rmk} 
 
We shall say that two self-energy clusters $T_{1},T_{2}$ have the same 
\emph{structure} if forcing $\nn_{\ell'_{T_{1}}}=\nn_{\ell'_{T_{2}}}$ one has 
$T_{1}=T_{2}$. Of course this provides an equivalence relation on the 
set of all self-energy clusters. The same consideration apply for left-fake 
and right-fake clusters. From now on we shall call 
self-energy, left-fake and right-fake clusters tout court such equivalence 
classes. 
 
A \emph{renormalised tree} is a tree in which no self-energy 
clusters appear; analogously a \emph{renormalised subgraph} is a 
subgraph of a tree $\theta$ which does not contains any self-energy cluster. 
Denote by $\Theta^{\RR}_{k,\nn}$ the set of renormalised trees 
with order $k$ and total momentum $\nn$, by 
$\gotR_{n}$ the set of renormalised self-energy clusters on scale $n$, 
and by $\gotL\gotF_{n}$ and $\gotR\gotF_{n}$ the sets of 
(renormalised) left-fake and right-fake clusters on scale $n$ respectively. 
 
For any $\theta\in\Theta^{\RR}_{k,\nn}$ we associate 
with each node $v\in N(\theta)$ a \emph{node factor} 
\begin{equation}\label{eq:3.9} 
\calF_{v}(\be_{0}) := \frac{1}{s_{v}!}\partial_{\be}^{s_{v}}F_{\nn_{v}}(\be_{0}). 
\end{equation} 
We associate with each line $\ell\in L(\theta)$ with $n_{\ell}\ge0$, a 
\emph{dressed propagator} $\calG_{n_{\ell}}(\oo\cdot\nn_{\ell};\e,\be_{0})$ 
(propagator tout court in the following) defined recursively as follows. 
 
\begin{figure}[ht] 
\centering 
\ins{360pt}{-132pt}{$x$} 
\ins{310pt}{-130pt}{$\displaystyle{\frac{\al_{0}}{4}}$} 
\ins{210pt}{-130pt}{$\displaystyle{\frac{\al_{0}}{8}}$} 
\ins{185pt}{-130pt}{$\displaystyle{\frac{\al_{m_{1}}}{4}}$} 
\ins{140pt}{-130pt}{$\displaystyle{\frac{\al_{m_{1}}}{8}}$} 
\ins{100pt}{-130pt}{$\displaystyle{\frac{\al_{m_{2}}}{8}}$} 
\ins{125pt}{-10pt}{$\Psi_{2}(x)$} 
\ins{200pt}{-10pt}{$\Psi_{1}(x)$} 
\ins{320pt}{-10pt}{$\Psi_{0}(x)$} 
\includegraphics[width=4in]{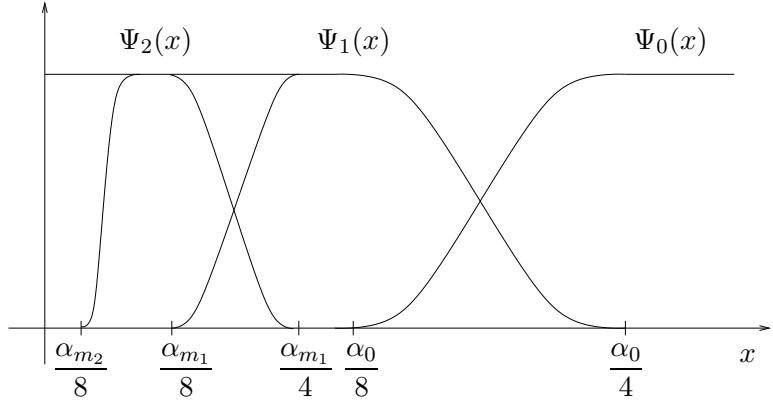} 
\vskip.2truecm 
\caption{Graphs of some of the $C^{\io}$
functions $\Psi_{n}(x)$ partitioning the unity in $\RRR\setminus\{0\}$; here 
$\al_{m}=\al_{m}(\oo)$. 
The function $\chi_{0}(x)=\chi(4x/\al_{0})$ is given by the sum of 
all functions 
$\Psi_{n}(x)$ for $n\ge 1$.} 
\label{fig:1} 
\end{figure}

Let us introduce the sequences $\{m_{n},p_{n}\}_{n \ge 0}$, with $m_{0}=0$ 
and, for all $n\ge 0$, $m_{n+1}=m_{n}+p_{n}+1$, 
where
$p_{n}:=\max\{q\in\ZZZ_{+}\,:\,\al_{m_{n}}(\oo)<2\al_{m_{n}+q}(\oo)\}$. Then the 
subsequence $\{\al_{m_{n}}(\oo)\}_{n\ge 0}$ of $\{\al_{m}(\oo)\}_{m\ge0}$ is 
decreasing. 
Let $\chi$ be a $C^{\io}$ non-increasing function such that 
\begin{equation}\label{eq:3.10} 
\chi(x)=\left\{ 
\begin{aligned} 
&1,\qquad |x| \le 1/2, \\ 
&0,\qquad |x| \ge 1. 
\end{aligned}\right. 
\end{equation} 
Set $\chi_{-1}(x)=1$ and $\chi_{n}(x)=\chi(4x/\al_{m_{n}}(\oo))$ for $n\ge0$. 
Set also $\psi(x)=1-\chi(x)$, $\psi_{n}(x)=\psi(4x/\al_{m_{n}}(\oo))$, 
and $\Psi_{n}(x)=\chi_{n-1}(x)\psi_{n}(x)$, for $n\ge 0$; see Figure 
\ref{fig:1}. 
 
\begin{lemma}\label{lem:3.1bis} 
For all $x\neq 0$ and for all $p\ge0$ one has 
\begin{equation} \nonumber 
\psi_{p}(x)+\sum_{n\ge p+1}\Psi_{n}(x)=1. 
\end{equation} 
\end{lemma} 
 
\prova For fixed $x\neq 0$ let $N=N(x):=\min\{n\,:\,\chi_{n}(x)= 0\}$ and 
note that $\max\{n\,:\,\psi_{n}(x)=0\} \le N-1$. Then if $p\le N-1$ 
\begin{equation} \nonumber 
\psi_{p}(x)+\sum_{n\ge p+1}\Psi_{n}(x)= \psi_{N-1}(x)+\chi_{N-1}(x)=1, 
\end{equation} 
while if $p\ge N$ one has 
\begin{equation} \nonumber 
\psi_{p}(x)+\sum_{n\ge p+1}\Psi_{n}(x)= \psi_{p}(x)=1. 
\end{equation} 
\EP 
\begin{rmk}\label{rmk:3.1bis} 
\emph{ 
Lemma \ref{lem:3.1bis} implies $\sum_{n\ge 0}\Psi_{n}(x)=1$ 
for all $x\neq0$. 
Hence $\{\Psi_{n}\}_{n\ge 0}$ is a partition of unity in $\RRR\setminus\{0\}$. 
} 
\end{rmk} 
 
Define, for $n\ge0$, 
\begin{equation}\label{eq:3.11} 
\calG_{n}(x;\e,\be_{0}) := \Psi_{n}(x) \left( x^{2}-\MM_{n-1}(x;\e,\be_{0})
\right)^{-1}, 
\end{equation} 
with formally, 
\begin{equation}\label{eq:3.13} 
\MM_{n-1}(x;\e,\be_{0}) := \sum_{q=-1}^{n-1}\chi_{q}(x)M_{q}(x;\e,\be_{0}), \qquad
M_{q}(x;\e,\be_{0}) := \sum_{T\in\gotR_{q}}\e^{k(T)}\Val_{T}(x;\e,\be_{0}), 
\end{equation} 
where $\Val_{T}(x;\e,\be_{0})$ is the \emph{renormalised value} 
of $T$, 
\begin{equation}\label{eq:3.14} 
\Val_{T}(x;\e,\be_{0}) := \left(\prod_{v\in N(T)}\calF_{v}(\be_{0}) 
\right) \left(\prod_{\ell\in L(T)}\calG_{n_{\ell}} 
(\oo\cdot\nn_{\ell};\e,\be_{0})\right). 
\end{equation} 
Here and henceforth, the sums and the products over empty sets have to be
considered as zero and $1$ respectively.
Note that $\Val_{T}$ depends on $\e$ because the
propagators do, 
and on $x=\oo\cdot\nn_{\ell'_{T}}$ only through the propagators 
associated with the lines $\ell\in\calP_{T}$ (see Remark \ref{rmk:3.1}).

\begin{rmk}\label{rmk:3.5bis}
\emph{
One has $|\gotR_{-1}|=1$ , so that
$\MM_{-1}(x;\e,\be_{0})=M_{-1}(x;\e,\be_{0})=
\e\partial_{\be_{0}}F_{\vzero}(\be_{0})$.
}
\end{rmk}
 
Set $\MM=\{\MM_{n}(x;\e,\be_{0})\}_{n\ge-1}$. 
We call \emph{self-energies} the quantities $\MM_{n}(x;\e,\be_{0})$. 
 
\begin{rmk}\label{rmk:3.1ter} 
\emph{ 
One has 
$$
\partial_{\be_{0}}\calG_{n}(x;\e,\be_{0})=\calG_{n}(x;\e,\be_{0}) 
\left( x^2-\MM_{n-1}(x;\e,\be_{0})\right)^{-1}
\partial_{\be_{0}}\MM_{n-1}(x;\e,\be_{0}).
$$ 
} 
\end{rmk} 
 
Set also $\calG_{-1}(0;\e,\be_{0})=1$ (so that we can associate 
a propagator also with the root line of $\theta\in 
\Theta_{k,\vzero}^{\RR}$). For any subgraph $S$ of any  
$\theta\in \Theta^{\RR}_{k,\nn}$ define the \emph{renormalised value} of $S$ as 
\begin{equation}\label{eq:3.15} 
\Val(S;\e,\be_{0}) := \left(\prod_{v\in N(S)}\calF_{v}(\be_{0}) 
\right) \left(\prod_{\ell\in L(S)}\calG_{n_{\ell}}(\oo\cdot\nn_{\ell}; 
\e,\be_{0})\right). 
\end{equation} 
Finally set 
\begin{equation}\label{eq:3.16} 
b_{\nn}^{[k]}(\e,\be_{0}) := \sum_{\theta\in \Theta^{\RR}_{k,\nn}}\Val 
(\theta;\e,\be_{0}), 
\end{equation} 
and 
\begin{equation}\label{eq:3.17} 
G^{[k]}(\e,\be_{0}) := \sum_{\theta\in \Theta^{\RR}_{k+1,\vzero}}\Val(\theta; 
\e,\be_{0}), 
\end{equation} 
and define formally 
\begin{equation}\label{eq:3.18} 
b^{\RR}(t;\e,\be_{0}) := \sum_{k\ge1}\e^{k}\sum_{\nn\in\ZZZ_{*}^{d}} 
e^{\ii\nn\cdot\oo t} b_{\nn}^{[k]}(\e,\be_{0}), 
\end{equation} 
and 
\begin{equation}\label{eq:3.19} 
G^{\RR}(\e,\be_{0}) := \sum_{k\ge 0}\e^{k} G^{[k]}(\e,\be_{0}). 
\end{equation} 
The series (\ref{eq:3.18}) and (\ref{eq:3.19}) will be called the 
\emph{resummed series}. 
The term ``resummed'' comes from the fact that if we formally 
expand (\ref{eq:3.18}) and (\ref{eq:3.19}) in powers of $\e$, we obtain 
(\ref{eq:3.1}) and (\ref{eq:3.4}) respectively, as it is easy to check. 
 
\begin{rmk}\label{rmk:3.2} 
\emph{ 
If $T$ is a renormalised left-fake (respectively right-fake) cluster, 
we can (and shall) write $\Val(T;\e,\be_{0})=\Val_{T}(\oo\cdot\nn_{\ell'_{T}}
;\e,\be_{0})$ since the propagators of the lines in $\calP_{T}$ depend on 
$\oo\cdot\nn_{\ell'_{T}}$. In particular one has 
$$ 
\sum_{T\in \gotL\gotF_{n}}\e^{k(T)}\Val_{T}(x;\e,\be_{0})= 
\sum_{T\in \gotR\gotF_{n}}\e^{k(T)}\Val_{T}(x;\e,\be_{0})= 
M_{n}(x;\e,\be_{0}). 
$$ 
} 
\end{rmk} 
 
\begin{rmk} \label{rmk:3.3} 
\emph{ 
Given a renormalised tree $\theta$ such that $\Val(\theta;\e,\be_{0}) 
\neq 0$, for any line $\ell\in L(\theta)$ (except possibly the root 
line) one has $\Psi_{n_{\ell}}(\oo\cdot\nn_{\ell})\neq0$, and hence 
\begin{equation} \label{eq:3.20} 
\frac{\al_{m_{n_{\ell}}}(\oo)}{8} < |\oo\cdot\nn_{\ell}| < 
\frac{\al_{m_{n_{\ell}-1}}(\oo)}{4} , 
\end{equation} 
where $\al_{m_{-1}}(\oo)$ has to be interpreted as $+\io$. 
The same considerations apply to any subgraph of $\theta$ and to 
renormalised self-energy clusters. 
Moreover, by the definition of $\{\al_{m_{n}}(\oo)\}_{n\ge 0}$, the 
number of scales which can be associated with a line $\ell$ in such a way 
that the propagator does not vanishes is at most 2; see Figure \ref{fig:1}. 
} 
\end{rmk} 
 
For $\theta\in\Theta^{\RR}_{k,\nn}$, let $\gotN_{n}(\theta)$ be the number 
of lines on scale $\ge n$ in $\theta$, and set 
\begin{equation}\label{eq:3.21} 
K(\theta):=\sum_{v\in N(\theta)}|\nn_{v}|. 
\end{equation} 
More generally, for any renormalised subgraph $T$ of any tree $\theta$ 
call $\gotN_{n}(T)$ the number of lines on scale $\ge n$ in $T$, and set 
\begin{equation}\label{eq:3.22} 
K(T):=\sum_{v\in N(T)}|\nn_{v}|. 
\end{equation} 
%
 
\begin{lemma}\label{lem:3.4} 
For any $\theta\in\Theta^{\RR}_{k,\nn}$ 
such that $\Val(\theta;\e,\be_{0}) \neq 0$ one has $\gotN_{n}(\theta) 
\le 2^{-(m_{n}-2)}K(\theta)$, for all $n\ge0$. 
\end{lemma} 
 
\prova 
First of all we note that if $\gotN_{n}(\theta)\ge 1$, then 
$K(\theta)\ge 2^{m_{n}-1}$. Indeed, if a line $\ell$ has scale $n_{\ell}\ge n$, 
then 
\begin{equation*} 
|\oo\cdot\nn_{\ell}|\le \frac{1}{4}\al_{m_{n-1}}(\oo)< 
\frac{1}{2}\al_{m_{n-1}+p_{n-1}}(\oo)=\frac{1}{2}\al_{m_{n}-1}(\oo)< 
\al_{m_{n}-1}(\oo), 
\end{equation*} 
and hence, by definition of $\al_{m}(\oo)$, one has 
$K(\theta)\ge |\nn_{\ell}|\ge 2^{m_{n}-1}$. 
Now we prove the bound 
$\gotN_{n}(\theta)\le \max\{2^{-(m_{n}-2)}K(\theta)-1,0\}$ by induction on 
the order. 
 
\begin{figure}[ht] 
\centering 
\ins{147pt}{-70pt}{$\ell_{\theta}$} 
\ins{142pt}{-50pt}{$\ge n$} 
\ins{214pt}{-058pt}{$< n $} 
\ins{250pt}{-17pt}{$\ge n$} 
\ins{270pt}{-28pt}{$\ell_{1}$} 
\ins{312pt}{006pt}{$\theta_{1}$} 
\ins{290pt}{-31pt}{$\ge n$} 
\ins{300pt}{-47pt}{$\ell_{2}$} 
\ins{330pt}{-10pt}{$\theta_{2}$} 
\ins{276pt}{-92pt}{$\ge n$} 
\ins{270pt}{-105pt}{$\ell_{r}$} 
\ins{320pt}{-104pt}{$\theta_{r}$} 
\includegraphics[width=3.0in]{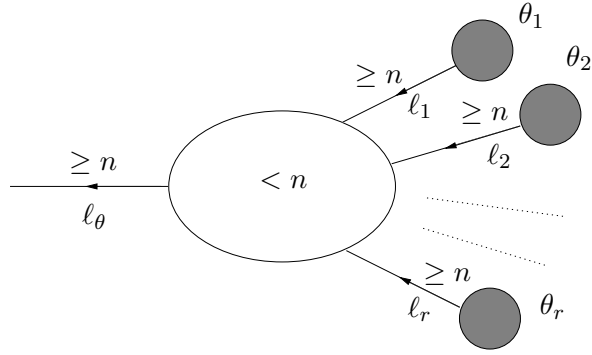} 
\caption{Construction used in the proof of Lemma \ref{lem:3.4} when 
$n_{\ell_{\theta}}\ge n$.} 
\label{fig:2} 
\end{figure} 
 
If the root line of $\theta$ has scale $n_{\ell_{\theta}}<n$ then the bound 
follows by the inductive hypothesis. If $n_{\ell_{\theta}}\ge n$, 
call $\ell_{1},\ldots,\ell_{r}$ the lines 
with scale $\ge n$ closest to $\ell_{\theta}$ (that is such that 
$n_{\ell'} < n$ for all lines $\ell'\in\calP(\ell_{\theta},\ell_{i})$, 
$i=1,\ldots,r$); see Figure \ref{fig:2}. If $r=0$ then $\gotN_{n}(\theta)=1$ 
and $|\nn| \ge 2^{m_{n}-1}$, so that the bound follows. 
If $r\ge 2$ the 
bound follows once more by the inductive hypothesis. If $r=1$, then 
$\ell_{1}$ is the only entering line of a cluster $T$ which is not a 
self-energy cluster as $\theta\in\Theta^{\RR}_{k,\nn}$, and hence 
$\nn_{\ell_{1}}\neq\nn$. But then 
$$ 
|\oo\cdot(\nn-\nn_{\ell_{1}})|\le|\oo\cdot\nn|+ 
|\oo\cdot\nn_{\ell_{1}}|\le \frac{1}{2}\al_{m_{n-1}}(\oo)< 
\al_{m_{n-1}+p_{n-1}}(\oo)=\al_{m_{n}-1}(\oo) , 
$$ 
as both $\ell_{\theta}$ and $\ell_{1}$ are on scale $\ge n$, so that one 
has $K(T)\ge|\nn-\nn_{\ell_{1}}|\ge 2^{m_{n}-1}$. Now, call 
$\theta_{1}$ the subtree of $\theta$ with root line $\ell_{1}$. Then one has 
\begin{equation*} 
\gotN_{n}(\theta)=1+\gotN_{n}(\theta_{1}) \le 1+ 
\max\{2^{-(m_{n}-2)}K(\theta_{1})-1,0\}, 
\end{equation*} 
so that 
\begin{equation*} 
\gotN_{n}(\theta) 
 \le 2^{-(m_{n}-2)}(K(\theta)-K(T))\le 2^{-(m_{n}-2)}K(\theta)-1, 
\end{equation*} 
again by induction. 
\EP 
 
\begin{lemma}\label{lem:3.5} 
For any $T\in\gotR_{n}$ such that $\Val_{T}(x;\e,\be_{0}) \neq 0$, 
one has $\gotN_{p}(T)\le 2^{-(m_{p}-2)}K(T)$, for all $0\le p\le n$. 
\end{lemma} 
 
\begin{figure}[ht] 
\centering 
\ins{228pt}{-38pt}{$\ell$} 
\ins{314pt}{-38pt}{$\ell_{T}'$} 
\ins{266pt}{-28pt}{$\widetilde{T}$} 
\ins{096pt}{-30pt}{$T =$} 
\includegraphics[width=3.0in]{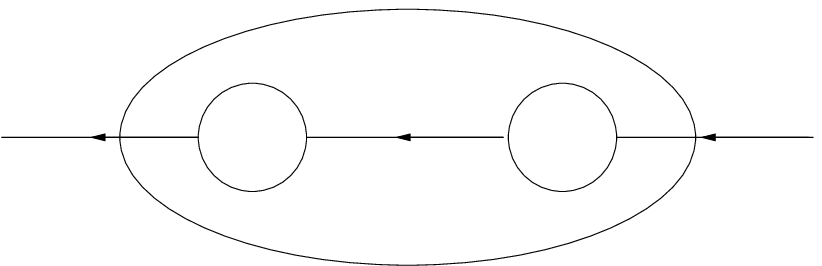} 
\caption{Construction used to prove 
$K(T) \ge 2^{m_{n}-1}$ when there is a line $\ell\in \calP_{T}$ on scale $n$.} 
\label{fig:3} 
\end{figure}

\prova 
We first prove that for all $n\ge 0$ and all $T\in\gotR_{n}$, one 
has $K(T)\ge 2^{m_{n}-1}$. In fact if $T\in\gotR_{n}$ then $T$ contains 
at least a line on scale $n$. If there is $\ell\in L(T)\setminus\calP_{T}$ 
with $n_{\ell}=n$, then 
$$ 
|\oo\cdot\nn_{\ell}|<\frac{1}{4}\al_{m_{n-1}}(\oo)<\al_{m_{n}-1}(\oo),
$$ 
and hence $K(T)\ge|\nn_{\ell}|>2^{m_{n}-1}$. Otherwise, let 
$\ell\in\calP_{T}$ be the line on scale $n$ which is closest to $\ell'_{T}$. 
Call $\widetilde{T}$ the subgraph (actually the cluster) consisting of all 
lines and nodes of $T$ preceding $\ell$; see Figure \ref{fig:3}. Then 
$\nn_{\ell}\neq \nn_{\ell'_{T}}$, otherwise $\widetilde{T}$ would be a 
self-energy cluster. Therefore $K(T)>|\nn_{\ell}-\nn_{\ell'_{T}}|>2^{m_{n}-1}$ 
as both $\ell,\ell'_{T}$ are on scale $\ge n$. 
 
\begin{figure}[ht]
\centering 
\ins{060pt}{-52pt}{$T =$} 
\ins{108pt}{-43pt}{$\ge p$} 
\ins{160pt}{-54pt}{$T_{1}$} 
\ins{220pt}{-43pt}{$\ge p$} 
\ins{220pt}{-63pt}{$\ell$} 
\ins{280pt}{-54pt}{$T_{2}$} 
\ins{334pt}{-43pt}{$\ge p$} 
\includegraphics[width=4.0in]{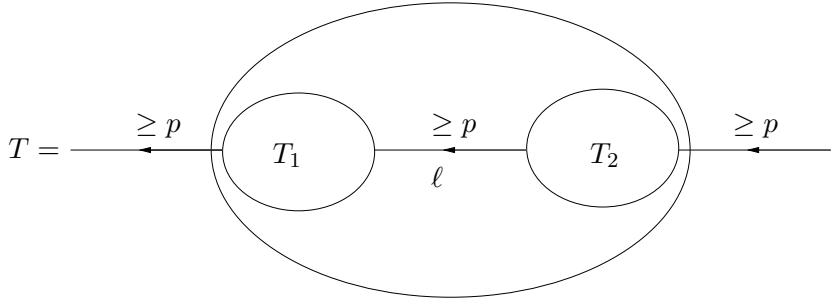} 
\caption{Construction used to prove Lemma \ref{lem:3.5}.} 
\label{fig:4} 
\end{figure}

Given a tree $\theta$, call $\CCCC(n,p)$ the set of  
renormalised subgraphs $T$ of 
$\theta$ with only one entering line $\ell'_{T}$ and one exiting line
$\ell_{T}$ both on scale $\ge p$, such that $L(T)\neq\emptyset$ and
$n_{\ell}\le n$ for any $\ell\in L(T)$. 
Note that $\gotR_{n}\subset \CCCC(n,p)$ for all $n,p\ge 0$ and,
reasoning as above,  one finds $K(T) \ge 2^{m_{q}-1}$, with $q=\min\{n,p\}$,
for all $T\in\CCCC(n,p)$. 
We prove that $\gotN_{p}(T)\le \max\{K(T)2^{-(m_{p}-2)}-1,0\}$ 
for all $0\le p\le n$ and all $T\in\CCCC(n,p)$. The proof is by induction on
the order. 
Call $N(\calP_{T})$ the set of nodes in $T$ connected by lines in 
$\calP_{T}$. If all lines in $\calP_{T}$ are on scale $< p$, then 
$\gotN_{p}(T)=\gotN_{p}(\theta_{1})+\ldots+\gotN_{p}(\theta_{r})$ 
if $\theta_{1},\ldots,\theta_{r}$ are the subtrees with root line 
entering a node in $N(\calP_{T})$, and hence the bound follows 
from (the proof of) Lemma \ref{lem:3.4}. 
If there exists a line $\ell\in\calP_{T}$ on scale $\ge p$, call 
$T_{1}$ and $T_{2}$ the subgraphs of $T$ such that 
$L(T)=\{\ell\}\cup L(T_{1})\cup L(T_{2})$, and note that if $L(T_{1}),L(T_{2})
\neq\emptyset$, then
$T_{1},T_{2}\in\CCCC(n,p)$; see Figure \ref{fig:4}. 
 
Hence, by the inductive hypothesis one has 
\begin{equation*} 
\begin{aligned} 
\gotN_{p}(T)&=1+\gotN_{p}(T_{1})+\gotN_{p}(T_{2})\\ 
&\le 1+\max\{2^{-(m_{p}-2)}K(T_{1})-1,0\}+\max\{2^{-(m_{p}-2)}K(T_{2})-1,0\}. 
\end{aligned} 
\end{equation*} 
If both $\gotN_{p}(T_{1}),\gotN_{p}(T_{2})$ are zero the bound trivially
follows as $K(T) \ge 2^{m_{p}-1}$, while 
if both are non-zero one has 
$$ 
\gotN_{p}(T)\le2^{-(m_{p}-2)}(K(T_{1})+K(T_{2}))-1 = 2^{-(m_{p}-2)}K(T)-1. 
$$ 
Finally if only one is zero, say $\gotN_{p}(T_{1})\ne0$ and
$\gotN_{p}(T_{2})=0$, 
$$ 
\gotN_{p}(T) \le 2^{-(m_{p}-2)}K(T_{1}) 
= 2^{-(m_{p}-2)}K(T)-2^{-(m_{p}-2)}K(T_{2}). 
$$ 
On the other hand, either $T_{2}\in \CCCC(n,p)$ or it is constituted by
only one node $v$ with $\nn_{v}\neq\vzero$, so that $K(T_{2})>2^{m_{p}-1}$ in
both cases. The same argument can be used in the case $\gotN_{p}(T_{1})=0$ 
and $\gotN_{p}(T_{2}) \neq 0$. 
\EP 
 
\zerarcounters 
\section{Convergence of the resummed series: part 1} 
\label{sec:4} 
 
To prove that the resummed series (\ref{eq:3.18}) converges, we first 
make the assumption that the propagators $\calG_{n_{\ell}} 
(x;\e,\be_{0})$ are bounded essentially as $1/x^{2}$: we shall see that 
in that case the convergence of the series can be easily proved. 
Then, in Section \ref{sec:5}, we shall check that the assumption is 
justified. 
 
\begin{defi}\label{defi:4.1} 
We shall say that $\MM$ satisfies \emph{property 1} if one has 
\begin{equation} \nonumber 
\Psi_{n+1}(x)|x^{2}-\MM_{n}(x;\e,\be_{0})|\ge \Psi_{n+1}(x)x^{2}/2, 
\end{equation} 
for all $n\ge -1$. 
\end{defi} 
 
\begin{lemma}\label{lem:4.2} 
Assume $\MM$ to satisfy property 1. Then 
the series (\ref{eq:3.18}) and (\ref{eq:3.19}) with the coefficients 
given by (\ref{eq:3.16}) and (\ref{eq:3.17}) respectively, converge 
for $\e$ small enough. 
\end{lemma} 
 
\prova 
Let $\theta\in\Theta^{\RR}_{k,\nn}$. 
The analyticity of $f$, hence of $F$, implies that there exist 
positive constants $F_{1},F_{2},\x$ such that for all $v\in N(\theta)$ one has 
\begin{equation} \label{eq:4.1} 
|\calF_{v}(\be_{0})|=\frac{1}{s_{v}!}|\partial_{\be}^{s_{v}} 
F_{\nn_{v}}(\be_{0})| 
\le F_{1}F_{2}^{s_{v}}e^{-\x|\nn_{v}|}. 
\end{equation} 
Moreover property 1 implies $|\calG_{n}(x;\e,\be_{0})|\le 2^{7} 
{\al_{m_{n}}(\oo)}^{-2}$ for all $n\ge 0$, and 
hence by Lemma \ref{lem:3.4} one can bound 
\begin{equation*} 
\begin{aligned} 
\prod_{\ell\in L(\theta)}|\calG_{n_{\ell}}(\oo\cdot\nn_{\ell};\e,\be_{0})| &\le 
 \prod_{n\ge 0}\left(\frac{2^{7}}{\al_{m_{n}}^{2}(\oo)}\right)^{N_{n}(\theta)} 
 \le\left(\frac{2^{7}}{\al_{m_{n_{0}}}^{2}(\oo)}\right)^{k} 
\prod_{n\ge n_{0}+1}\left(\frac{2^{7}}{\al_{m_{n}}^{2}(\oo)} 
\right)^{N_{n}(\theta)}\\ 
 &\le \left(\frac{2^{7}}{\al_{m_{n_{0}}}^{2}(\oo)}\right)^{k} 
\prod_{n\ge n_{0}+1}\left(\frac{2^{7/2}}{\al_{m_{n}}(\oo)}\right)^{2^{-(m_{n}-3)} 
K(\theta)}\\ 
 &\le \left(\frac{2^{7}}{\al_{m_{n_{0}}}^{2}(\oo)}\right)^{k} 
{\rm exp}\left(8K(\theta) 
\sum_{n\ge n_{0}+1}\frac{1}{2^{m_{n}}} 
\log\frac{2^{7/2}}{\al_{m_{n}}(\oo)}\right) \\ 
 &\le D^{k}(n_{0}){\rm exp}(\x(n_{0})K(\theta)), 
\end{aligned} 
\end{equation*} 
with 
$$ 
D(n_{0})=\frac{2^{7}}{\al_{m_{n_{0}}}^{2}(\oo)},\qquad 
\x(n_{0})=8\sum_{n\ge n_{0}+1}\frac{1}{2^{m_{n}}}\log\frac{2^{7/2}}{\al_{m_{n}} 
(\oo)}. 
$$ 
Then, by Hypothesis \ref{hyp1}, one can choose $n_{0}$ such that 
$\x(n_{0})\le \x/2$. The sum over the other labels is bounded by a 
constant to the power $k$, and hence 
one can bound 
$$ 
\sum_{\theta\in\Theta^{\RR}_{k,\nn}}|\Val(\theta;\e,\be_{0})|\le 
C_{0}C_{1}^{k}e^{-\x |\nn|/2}, 
$$ 
for some constants $C_{0}, C_{1}$, and this is enough to prove the 
assertion. 
\EP 
 
\begin{lemma}\label{lem:4.3} 
Assume $\MM$ to satisfy property 1. Then for $\e$ small enough 
the function (\ref{eq:3.18}), 
with the coefficients given by (\ref{eq:3.16}), solves the equation 
(\ref{eq:2.3a}). 
\end{lemma} 
 
\prova 
We shall prove that, the function $b^{\RR}$ defined in (\ref{eq:3.18}) 
satisfies the equation of motion (\ref{eq:2.3a}), i.e. we shall check that 
$b^{\RR}=\e g F(\oo t,\be_{0}+b^{\RR})$, where $g$ is the 
pseudo-differential operator with kernel $g(\oo\cdot\nn)= 
1/(\oo\cdot\nn)^{2}$. We can write the Fourier coefficients of $b^{\RR}$ 
as 
\begin{equation} \label{eq:4.2} 
b^{\RR}_{\nn}=\sum_{n\ge 0}b_{\nn}^{[n]},\qquad 
b_{\nn}^{[n]}=\sum_{k\ge1}\e^{k}\sum_{\theta \in \Theta^{\RR}_{k,\nn}(n)} 
\Val(\theta;\e,\be_{0}), 
\end{equation} 
where $\Theta^{\RR}_{k,\nn}(n)$ is the subset of $\Theta^{\RR}_{k,\nn}$ 
such that $n_{\ell_{\theta}}=n$. 
 
Using Remark \ref{rmk:3.1bis} and Lemma \ref{lem:4.2}, in Fourier space one 
can write 
\begin{equation} \nonumber 
\begin{aligned} 
g(\oo\cdot\nn)&[\e F(\oo t,\be_{0}+b^{\RR})]_{\nn} 
= g(\oo\cdot\nn)\sum_{n\ge0}\Psi_{n}(\oo\cdot\nn) 
 [\e F(\oo t,\be_{0}+b^{\RR})]_{\nn} \\ 
&= g(\oo\cdot\nn)\sum_{n\ge0}\Psi_{n}(\oo\cdot\nn) 
(\calG_{n}(\oo\cdot\nn;\e,\be_{0}))^{-1}\calG_{n}(\oo\cdot\nn;\e,\be_{0}) 
 [\e F(\oo t,\be_{0}+b^{\RR})]_{\nn} \\ 
&= g(\oo\cdot\nn)\sum_{n\ge0}\left((\oo\cdot\nn)^{2}-\MM_{n-1}(\oo\cdot\nn; 
 \e,\be_{0})\right)\calG_{n}(\oo\cdot\nn;\e,\be_{0}) 
 [\e F(\oo t,\be_{0}+b^{\RR})]_{\nn} \\ 
&= g(\oo\cdot\nn)\sum_{n\ge0}\left((\oo\cdot\nn)^{2}-\MM_{n-1}(\oo\cdot\nn; 
 \e,\be_{0})\right)\sum_{k\ge1}\e^{k}\sum_{\theta\in\ol{\Theta}^{\RR}_{k,\nn}(n)} 
 \Val(\theta;\e,\be_{0}), 
\end{aligned} 
\end{equation} 
where $\ol{\Theta}^{\RR}_{k,\nn}(n)$ differs from $\Theta^{\RR}_{k,\nn}(n)$ 
as it contains also trees $\theta$ which have one self-energy 
cluster with exiting line $\ell_{\theta}$. 
If we separate the trees containing such self-energy cluster from 
the others, we obtain 
\begin{equation} \nonumber 
\begin{aligned} 
g(\oo\cdot\nn)&[\e F(\oo t,\be_{0}+b^{\RR})]_{\nn} 
= g(\oo\cdot\nn)\sum_{n\ge0}\left((\oo\cdot\nn)^{2}-\MM_{n-1}(\oo\cdot\nn; 
 \e,\be_{0})\right)b_{\nn}^{[n]} \\ 
&\quad+g(\oo\cdot\nn)\sum_{n\ge0}\Psi_{n}(\oo\cdot\nn) 
 \sum_{p\ge n}\sum_{q=-1}^{n-1}M_{q}(\oo\cdot\nn;\e,\be_{0})b_{\nn}^{[p]} \\ 
&\quad+g(\oo\cdot\nn)\sum_{n\ge1}\Psi_{n}(\oo\cdot\nn) 
 \sum_{p=0}^{n-1}\sum_{q=-1}^{p-1}M_{q}(\oo\cdot\nn;\e,\be_{0}) 
 b_{\nn}^{[p]} \\ 
&= g(\oo\cdot\nn)\sum_{n\ge0}\left((\oo\cdot\nn)^{2}-\MM_{n-1}(\oo\cdot\nn; 
 \e,\be_{0})\right)b_{\nn}^{[n]} \\ 
&\quad+g(\oo\cdot\nn)\sum_{p\ge0}\left(\sum_{q=-1}^{p-1} 
 M_{q}(\oo\cdot\nn;\e,\be_{0})\sum_{n\ge q+1}\Psi_{n}(\oo\cdot\nn) 
 \right)b_{\nn}^{[p]} \\ 
&=  g(\oo\cdot\nn)\sum_{n\ge0}\left((\oo\cdot\nn)^{2}-\MM_{n-1}(\oo\cdot\nn; 
 \e,\be_{0})\right)b_{\nn}^{[n]} \\ 
&\quad+g(\oo\cdot\nn)\sum_{n\ge0}\left(\sum_{q=-1}^{n-1} 
 M_{q}(\oo\cdot\nn;\e,\be_{0})\chi_{q}(\oo\cdot\nn)\right)b_{\nn}^{[n]} \\ 
&=  g(\oo\cdot\nn)\sum_{n\ge0}\left((\oo\cdot\nn)^{2}-\MM_{n-1}(\oo\cdot\nn; 
 \e,\be_{0})\right)b_{\nn}^{[n]} \\ 
&\quad+g(\oo\cdot\nn)\sum_{n\ge0}\MM_{n-1}(\oo\cdot\nn;\e,\be_{0}) 
 b_{\nn}^{[n]} \\ 
&=\sum_{n\ge0}b_{\nn}^{[n]}=b^{\RR}_{\nn}, 
\end{aligned} 
\end{equation} 
so that the proof is complete. 
\EP 
 
\begin{defi}\label{defi:4.4} 
We shall say that $\MM$ satisfies \emph{property 2-$p$} if one has 
\begin{equation} \nonumber 
\Psi_{n+1}(x)|x^{2}-\MM_{n}(x;\e,\be_{0})|\ge \Psi_{n+1}(x)x^{2}/2, 
\end{equation} 
for all $-1\le n <p$. 
\end{defi} 
 
\begin{lemma}\label{lem:4.6} 
Assume $\MM$ to satisfy property 2-$p$. Then for any 
$0\le n\le p$ the self-energies are well defined and one has 
\begin{subequations} 
\begin{align} 
& |M_{n}(x;\e,\be_{0})|\le \e^{2}K_{1}e^{-K_{2}2^{m_{n}}}, 
\label{eq:4.3a} \\ 
& |\partial^{j}_{x}M_{n}(x;\e,\be_{0})|\le 
\e^{2}C_{j}e^{-\ol{C}_{j}2^{m_{n}}}, 
\qquad j=1,2, 
\label{eq:4.3b} 
\end{align} 
\label{eq:4.3} 
\end{subequations} 
for suitable constants $K_{1},K_{2},C_{1},C_{2},\ol{C}_{1}$ and 
$\ol{C}_{2}$. 
\end{lemma} 
 
\prova 
Property 2-$p$ implies 
$|\calG_{n}(x;\e,\be_{0})|\le 2^{7} \al_{m_{n}}(\oo)^{-2}$ 
for all $0\le n\le p$. Then, using also Lemma \ref{lem:3.5} 
and the fact that any self-energy cluster in $\gotR_{n}$ has at least 
two nodes for any $n\ge0$, we obtain 
$$ 
|M_{n}(x;\e,\be_{0})|\le \sum_{T\in\gotR_{n}}|\e|^{k(T)}|\Val_{T}(x;\e,\be_{0})| 
\le\sum_{k\ge 2}|\e|^{k}C^{k}e^{-K_{2}2^{m_{n}}}, 
$$ 
so that (\ref{eq:4.3a}) is proved for $\e$ small enough. 
Now we prove (\ref{eq:4.3b}) by induction on $n$. 
For $n=0$ the bound is obvious. Assume then (\ref{eq:4.3b}) to hold for 
all $n'<n$. For any $T\in\gotR_{n}$ such that $\Val_{T}(x;\e,\be_{0})\ne0$ 
one has 
\begin{equation} \nonumber 
\partial_{x}\Val_{T}(x;\e,\be_{0})= 
\sum_{\ell\in\calP_{T}}\!\!\left(\prod_{v\in N(T)} 
\calF_{v}(\be_{0})\!\!\right)\!\! 
\left(\partial_{x}\calG_{n_{\ell}}(x_{\ell}; 
\e,\be_{0})\!\!\!\!\!\! 
\prod_{\ell'\in L(T)\setminus\{\ell\}} 
\!\!\!\!\!\! 
\calG_{n_{\ell'}}(\oo\cdot\nn_{\ell'};\e,\be_{0})\!\!\right), 
\end{equation} 
where $x_{\ell}=\oo\cdot\nn_{\ell}=x+\oo\cdot\nn_{\ell}^{0}$ and 
\begin{equation} \nonumber 
\begin{aligned} 
\partial_{x}\calG_{n_{\ell}}&(x_{\ell};\e,\be_{0})= 
\frac{d}{dx}\calG_{n_{\ell}} 
(\oo\cdot\nn_{\ell}^{0}+x;\e,\be_{0})\\ 
&=\frac{\partial_{x}\Psi_{n_{\ell}}(x_{\ell})} 
{x_{\ell}^{2}-\MM_{n_{\ell}-1}(x_{\ell};\e,\be_{0})}- 
\frac{\Psi_{n_{\ell}}(x_{\ell}) 
\left(2x_{\ell}-\partial_{x}\MM_{n_{\ell}-1}(x_{\ell};\e,\be_{0}) 
\right)} 
{\left(x_{\ell}^{2}-\MM_{n_{\ell}-1}(x_{\ell};\e,\be_{0}) 
\right)^{2}}. 
\end{aligned} 
\end{equation} 
One has 
\begin{equation} \nonumber 
|\partial_{x}\Psi_{n_{\ell}}(x_{\ell})|\le |\partial_{x}\chi_{n_{\ell}-1} 
    (x_{\ell})|+|\partial_{x}\psi_{n_{\ell}}(x_{\ell})|\le 
    \frac{B_{1}}{\al_{m_{n_{\ell}}}(\oo)}, 
\end{equation} 
for some constant $B_{1}$
and, by (\ref{eq:4.3a}), the inductive hypothesis and Hypothesis \ref{hyp1}, 
\begin{equation} \nonumber 
\begin{aligned} 
|\partial_{x}\MM_{n_{\ell}-1}(x_{\ell};\e,\be_{0})| 
    &\le \sum_{q=0}^{n_{\ell}-1}|(\partial_{x}\chi_{q}(x_{\ell})) 
    M_{q}(x_{\ell};\e,\be_{0})| 
    +\sum_{q=0}^{n_{\ell}-1}|\partial_{x}M_{q}(x_{\ell};\e,\be_{0})|\\ 
 &\le\e^{2}B_{1}K_{1}\sum_{q\ge 0}\frac{1}{\al_{m_{q}}(\oo)} 
    e^{-K_{2}2^{m_{q}}}+\e^{2}C_{1}\sum_{q\ge 0}e^{-\ol{C}_{1}2^{m_{q}}}\\ 
 &\le \e^{2}B_{2}, 
\end{aligned} 
\end{equation} 
for some constant $B_{2}$. 
Hence, at the cost of replacing the bound for the propagators with 
$\widetilde{C}\al_{m_{n_{\ell}}}(\oo)^{-4}$ for some constant $\widetilde{C}$, 
one can rely upon Lemma \ref{lem:3.5} to obtain (\ref{eq:4.3b}) for $j=1$. 
For $j=2$ one can reason analogously. 
\EP 
 
\begin{lemma}\label{lem:4.5} 
Assume $\MM$ to satisfy property 2-$p$. Then one has 
$\MM_{n}(x;\e,\be_{0})=\MM_{n}(0;\e,\be_{0})+O(\e^{2}x^{2})$ 
for all $0 \le n\le p$. 
\end{lemma} 
 
\prova 
We shall prove that $\MM_{n}(x;\e,\be_{0})=\MM_{n}(-x;\e,\be_{0})$, 
by induction on $n\ge-1$. For $n=-1$ the identity is obvious 
since $\MM_{-1}$ does not depend on $x$. 
Assume now $\MM_{q}(x;\e,\be_{0})=\MM_{q}(-x;\e,\be_{0})$ for 
all $q<n$. This implies $\calG_{q}(x;\e,\be_{0})=\calG_{q}(-x;\e,\be_{0})$ 
for $q\le n$. Let $T\in\gotR_{n}$ and consider the self-energy cluster 
$T_{1}$ obtained from $T$ by taking $\ell_{T}$ as the entering line and 
$\ell'_{T}$ as the exiting line (i.e. $\ell_{T_{1}}'=\ell_{T}$ and 
$\ell_{T_{1}}=\ell_{T}'$) and by taking $\nn_{\ell'_{T_{1}}}= 
-\nn_{\ell'_{T}}$. Hence the momenta of the lines belonging to $\calP_{T}$ 
change signs, while all the other momenta do not change: therefore 
all propagators are left unchanged. Hence $M_{n}(x;\e,\be_{0})= 
M_{n}(-x;\e,\be_{0})$, so that 
$\partial_{x}M_{n}(0;\e,\be_{0})=0$ for all $n \le p$, and, by Lemma 
\ref{lem:4.6}, this is enough to prove the assertion. 
\EP 
 
\begin{lemma}\label{lem:4.7} 
Assume $\MM$ to satisfy property 1. 
Then the function $G^{\RR}(\e,\be_{0})$ and the self-energies 
$\MM_{n}(x;\e,\be_{0})$ are $C^{\io}$ in both $\e$ and $\be_{0}$. 
\end{lemma} 
 
\prova 
It follows from the explicit expressions for $G^{\RR}(\e,\be_{0})$ and 
$\MM_{n}(x;\e,\be_{0})$. 
\EP 
 
Define formally 
\begin{equation}\label{eq:4.4} 
\MM_{\io}(x;\e,\be_{0})=\lim_{n\to\io}\MM_{n}(x;\e,\be_{0}), 
\end{equation} 
and note that if $\MM$ satisfies property 1, then $\MM_{\io}(x;\e,\be_{0})$ 
is well defined and moreover it is $C^{\io}$ in both $\e$ and $\be_{0}$. 
 
The following result plays a crucial role. The proof is 
deferred to Appendix \ref{app:a}. 
 
\begin{lemma}\label{lem:4.8} 
Assume $\MM$ to satisfy property 1. Then one has 
\begin{equation}  \nonumber 
\e\partial_{\be_{0}} G^{\RR}(\e,\be_{0})=\MM_{\io}(0;\e,\be_{0}). 
\end{equation} 
\end{lemma} 
 
\begin{rmk} \label{rmk:4.9} 
\emph{ 
If we take the formal power expansions of both $G^{\RR}(\e,\be_{0})$ and 
$\MM_{\io}(0;\e,\be_{0})$, 
we obtain tree expansions where self-energy clusters are allowed; see 
Section \ref{sec:6} for further details. Then 
the identity $\e\partial_{\be_{0}} 
G^{\RR}(\e,\be_{0})=\MM_{\io}(0;\e,\be_{0})$ 
is easily found to be 
satisfied to any perturbation order. However, without any resummation 
procedure, we are no longer able to prove the convergence of the 
series, so that the identity becomes a meaningless ``$\infty=\infty$''. 
} 
\end{rmk} 
 
\begin{rmk} \label{rmk:4.10} 
\emph{ 
The identity $\e\partial_{\be_{0}}G^{\RR}(\e,\be_{0})=\MM_{\io}(0;\e,\be_{0})$, 
in Lemma \ref{lem:4.8},  can be seen as an identity between classes of
diagrams. In turn, in light of a possible quantum field formulation of the
problem, this can be thought as a consequence of some deep Ward identity
of the corresponding field theory. Ward identities
play a crucial role in quantum field theory. The analogy 
between KAM theory and quantum field theory has been widely stressed 
in the literature \cite{GGM,BGK,DK}; in particular the cancellations 
which assure the convergence of the perturbation series for maximal 
KAM tori are deeply related to a Ward identity, as shown in 
\cite{BGK}, which can be seen as a remarkable identity between classes 
of graphs. In the case studied in this paper, we have a similar 
situation, made fiddlier by the fact that we have to deal with 
nonconvergent series to be resummed, and it is well known that 
identities which are trivial on a formal level can turn out to be 
difficult to prove rigorously \cite{M}. However, we expect
a Ward identity to hold also in our case, so as to imply
that $\e\partial_{\be_{0}}G^{\RR}(\e,\be_{0})=\MM_{\io}(0;\e,\be_{0})$.
It would be interesting to confirm the expectation and to determine
the Ward identity explicitly.
} 
\end{rmk} 
 
\begin{lemma}\label{lem:4.11} 
Assume $\MM$ to satisfy property 1. Then the implicit function 
equation $G^{\RR}(\e,\be_{0})=0$ admits a solution $\be_{0}=\be_{0}(\e)$, 
such that $\be_{0}(0)=\be_{0}^{*}$. 
Moreover in a suitable half-neighbourhood of $\e=0$, one has 
$\e\partial_{\be_{0}}G^{\RR}(\e,\be_{0}(\e))\le 0$. 
\end{lemma} 
 
\prova 
Property 1 allows us to write $G^{\RR}(\e,\be_{0})=F_{\vzero}(\be_{0})+ 
O(\e)$, so that by Hypothesis \ref{hyp2} one has 
$\partial_{\be_{0}}^{\gotn}G^{\RR}(0,\be_{0}^{*})\neq 0$. Then there exist 
two half-neighbourhood $V_{-},V_{+}$ of 
$\be_{0}=\be_{0}^{*}$ such that $G^{\RR}(0,\be_{0})>0$ for 
$\be_{0}\in V_{+}$ and $G^{\RR}(0,\be_{0})<0$ for $\be_{0}\in V_{-}$. 
Hence, by continuity, for all $\be_{0}\in V_{+}$ there exists a 
neighbourhood $U_{+}(\be_{0})$ of $\e=0$ such that $G^{\RR}(\e,\be_{0})>0$ 
for all $\e\in U_{+}(\be_{0})$ 
and, for the same reason, for all $\be_{0}\in V_{-}$ there exists a 
neighbourhood $U_{-}(\be_{0})$ of $\e=0$ such that $G^{\RR}(\e,\be_{0})<0$ 
for all $\e\in U_{-}(\be_{0})$. Therefore, again by continuity, there 
exists a continuous curve $\be_{0}=\be_{0}(\e)$ defined in a suitable 
neighbourhood $U=(-\ol{\e},\ol{\e})$ such that 
$\be_{0}(0)=\be_{0}^{*}$ and $G^{\RR}(\e,\be_{0}(\e))\equiv0$. 
Moreover, if $\partial_{\be_{0}}^{\gotn}G^{\RR}(0,\be_{0}^{*})>0$, then 
$V_{+},V_{-}$ are of the form $(\be_{0}^{*},v_{+})$ and $(v_{-},\be_{0}^{*})$ 
respectively, and therefore 
$\partial_{\be_{0}}G^{\RR}(c,\be_{0}(c))\ge 0$ for all $c\in U$. 
If on the contrary $\partial_{\be_{0}}^{\gotn}G^{\RR}(0,\be_{0}^{*})<0$, 
one has $V_{+}=(v_{+},\be_{0}^{*})$ and $V_{-}=(\be_{0}^{*},v_{-})$, and then 
$\partial_{\be_{0}}G^{\RR}(c,\be_{0}(c))\le 0$ for all $c\in U$. 
Hence the assertion follows in both cases, again by Hypothesis \ref{hyp2}. 
\EP 
 
\begin{rmk}\label{rmk:4.12} 
\emph{ 
If $\MM$ satisfies property 1, one has 
$$ 
G^{\RR}(\e,\be_{0})=[F(\oo t,\be_{0}+b^{\RR}(t;\e,\be_{0}))]_{\vzero}, 
$$ 
and hence, if $\be_{0}=\be_{0}(\e)$ is the solution referred to in Lemma 
\ref{lem:4.11}, by Lemma \ref{lem:4.3} the function 
$$ 
\be(t;\e)=\be_{0}(\e)+b^{\RR}(t;\e,\be_{0}(\e)), 
$$ 
solves the equation of motion (\ref{eq:1.1}). 
} 
\end{rmk} 
 
\begin{rmk}\label{rmk:4.13} 
\emph{ 
In Lemma \ref{lem:4.11} we widely used that the variable $\be_{0}$ is
one-dimensional. All the other results in this paper could be quite
easily extended to higher dimension.
} 
\end{rmk} 
 
\begin{rmk}\label{rmk:4.14} 
\emph{The results of this section are not sufficient to prove Theorem
\ref{thm:2.1} because we have assumed -- without proving -- that property
1 is satisfied. In Section \ref{sec:5} we shall show that, thanks to the
symmetry property of Lemma \ref{lem:4.5} and the identity of Lemma 
\ref{lem:4.8}, property 1 is satisfied along a suitable continuous curve
$\be_{0}=\ol{\be}_{0}(\e)$ such that $G^{\RR}(\e,\ol{\be}_{0}(\e))=0$. 
} 
\end{rmk} 
 
\zerarcounters 
\section{Convergence of the resummed series: part 2} 
\label{sec:5} 
 
In this section we shall remove the assumption that the self-energies 
satisfy property 1 of Definition \ref{defi:4.1} -- see Remark \ref{rmk:4.14}. 
For all $n\ge 0$, define the $C^{\io}$ non-increasing functions $\x_{n}$ 
such that 
\begin{equation}\label{eq:5.1} 
\x_{n}(x)=\left\{ 
\begin{aligned} 
&1,\quad x\le \al_{m_{n+1}}^{2}(\oo)/2^{9},\\ 
&0,\quad x\ge \al_{m_{n+1}}^{2}(\oo)/2^{8}, 
\end{aligned}\right. 
\end{equation} 
and set $\x_{-1}(x)=1$. Define recursively, for all $n\ge0$, the propagators 
\begin{equation}\label{eq:5.2} 
\ol{\calG}_{n}(x;\e,\be_{0})=\Psi_{n}(x) 
\left( x^{2}-\ol{\MM}_{n-1}(x;\e,\be_{0})\x_{n-1}(\ol{\MM}_{n-1} 
(0;\e,\be_{0})) \right)^{-1}, 
\end{equation} 
with $\ol{\MM}_{-1}(x;\e,\be_{0})=\e\partial_{\be}F_{\vzero}(\be_{0})$, 
and for $n\ge0$ 
\begin{equation}\label{eq:5.3} 
\ol{\MM}_{n}(x;\e,\be_{0})=\ol{\MM}_{n-1}(x;\e,\be_{0})+ 
\chi_{n}(x)\ol{M}_{n}(x;\e,\be_{0}), 
\end{equation} 
where we have set 
\begin{equation}\label{eq:5.4} 
\ol{M}_{n}(x;\e,\be_{0})=\sum_{T\in\gotR_{n}}\e^{k(T)} 
\ol{\Val}_{T}(x;\e,\be_{0}), 
\end{equation} 
with 
\begin{equation}\label{eq:5.5} 
\ol{\Val}_{T}(x;\e,\be_{0})= 
\left(\prod_{v\in N(T)}\calF_{v}(\be_{0})\right)\left(\prod_{\ell\in L(T)} 
\ol{\calG}_{n_{\ell}}(\oo\cdot\nn_{\ell};\e,\be_{0})\right), 
\end{equation} 
and $x=\oo\cdot\nn_{\ell'_{T}}$. 
 
Set also 
$\ol{\MM}=\{\ol{\MM}_{n}(x;\e,\be_{0})\}_{n\ge-1}$, and 
$\ol{\MM}^{\x}=\{\ol{\MM}_{n}(x;\e,\be_{0})\x_{n}(\ol{\MM}_{n} 
(0;\e,\be_{0}))\}_{n\ge-1}$. 
 
\begin{lemma}\label{lem:5.1} 
$\ol{\MM}^{\x}$ satisfies property 1.
\end{lemma}

\prova
We shall prove that $\ol{\MM}^{\x}$ satisfies property 2-$p$ for all
$p\ge0$, by induction on $p$. Property 2-0 is trivially satisfied for
$\e$ small enough. Assume $\ol{\MM}^{\x}$ to satisfy
property 2-$p$. Then we can repeat (almost word by word) the proofs of
Lemmas \ref{lem:4.6} and \ref{lem:4.5} so as to obtain
$$
\ol{\MM}_{p}(x;\e,\be_{0})=\ol{\MM}_{p}(0;\e,\be_{0})+O(\e^{2}x^{2}),
$$
hence, by the definition of the function $\x_{p}$, $\ol{\MM}^{\x}$
satisfies property 2-$(p+1)$, and thence the assertion follows.
\EP
 
Set 
\begin{equation}\label{eq:5.6} 
\ol{\Val}(\theta;\e,\be_{0})=\left(\prod_{v\in N(\theta)}\calF_{v}(\be_{0}) 
\right) \left(\prod_{\ell\in L(\theta)}\ol{\calG}_{n_{\ell}}(\oo\cdot\nn_{\ell}; 
\e,\be_{0})\right), 
\end{equation} 
and 
\begin{equation}\label{eq:5.7} 
\ol{b}^{[k]}_{\nn}(\e,\be_{0})=\sum_{\theta\in\Theta^{\RR}_{k,\nn}}\ol{\Val} 
(\theta;\e,\be_{0}), 
\end{equation} 
and define 
\begin{equation}\label{eq:5.8} 
\ol{b}(t,\e,\be_{0})=\sum_{k\ge1}\e^{k} 
\ol{b}^{[k]}(\e,\be_{0})=\sum_{k\ge 1} \e^{k}\sum_{\nn\in\ZZZ^{d}_{*}} 
e^{\ii \nn\cdot\oo t}\ol{b}^{[k]}_{\nn}(\e,\be_{0}). 
\end{equation} 
Note that, by (the proof of) Lemma \ref{lem:4.2} the series (\ref{eq:5.8}) 
converges.

Define also
\begin{equation}\label{eq:5.9a}
\ol{\MM}_{\io}(x;\e,\be_{0}):=\lim_{n\to\io}
\ol{\MM}_{n}(x;\e,\be_{0}),
\end{equation}
and note that, by Lemma \ref{lem:5.1} the limit in (\ref{eq:5.9a}) is well
defined and it is $C^{\io}$ in both $\e$ and $\be_{0}$.
Introduce the $C^{\io}$ functions $\ol{G}(\e,\be_{0})$ such
that  $\ol{\MM}_{\io}(0;\e,\be_{0})=\e\partial_{\be_{0}}\ol{G}(\e,\be_{0})$
and $\ol{G}(0,\be_{0}^{*})=0$, and for any such function consider the
implicit function equation 
\begin{equation}\label{eq:5.9} 
\ol{G}(\e,\be_{0})=0. 
\end{equation} 
%
 
\begin{lemma}\label{lem:5.2} 
The implicit function equation (\ref{eq:5.9}) admits a solution 
${\be}_{0}=\ol{\be}_{0}(\e)$ such that $\ol{\be}_{0}(0)=\be_{0}^{*}$. 
Moreover in a suitable half-neighbourhood of $\e=0$, one has 
$\e\partial_{\be_{0}} \ol{G}(\e,\ol{\be}_{0}(\e))\le0$. 
\end{lemma} 
 
\prova 
By construction, all the functions $\ol{G}(\e,\be_{0})$ are smooth and
of the form $\ol{G}(\e,\be_{0})=F_{\vzero}(\be_{0})+O(\e)$. Then
the result follows straightforward from (the proof of) Lemma \ref{lem:4.11}. \EP 
 
\begin{lemma}\label{lem:5.3} 
Let $\be_{0}=\ol{\be}_{0}(\e)$ be the solution referred to in Lemma 
\ref{lem:5.2}. Then 
one has $\x_{n}(\ol{\MM}_{n}(0;\e,\ol{\be}_{0}(\e)))\equiv 1$ 
for all $n\ge0$, in a suitable half-neighbourhood of $\e=0$. 
\end{lemma} 
 
\prova 
If ${\be}_{0}=\ol{\be}_{0}(\e)$, by Lemma \ref{lem:5.2} in a suitable 
half-neighbourhood of $\e=0$ one has $\ol{\MM}_{\io}(0;\e,\ol{\be}_{0}(\e))= 
\e\partial_{\be_{0}}\ol{G}(\e,\ol{\be}_{0}(\e))\le 0$. Hence, 
as the bound (\ref{eq:4.3a}) holds also for $\ol{M}_{n}(x;\e,\be_{0})$, 
one has 
\begin{equation}\label{eq:5.10} 
\begin{aligned} 
\ol{\MM}_{n}(0;\e,\ol{\be}_{0}(\e)) & 
  \le\ol{\MM}_{n}(0;\e,\ol{\be}_{0}(\e))- 
  \ol{\MM}_{\io}(0;\e,\ol{\be}_{0}(\e))\\ 
&\le \sum_{p\ge n+1}|\ol{M}_{p}(0;\e,\ol{\be}_{0}(\e))| \\ 
&\le 2K_{1}\e^{2}e^{-K_{2}2^{m_{n}}} \le \frac{\al_{m_{n+1}}^{2}}{2^{11}}, 
\end{aligned} 
\end{equation} 
so that the assertion follows by the definition of $\x_{n}$. 
\EP 
 
\begin{lemma}\label{lem:5.4} 
For $\be_{0}=\ol{\be}_{0}(\e)$, one has $\MM=\ol{\MM}=\ol{\MM}^{\x}$, and 
hence one can choose $\ol{G}(\e,\be_{0})$ such that
$G^{\RR}(\e,\ol{\be}_{0}(\e))=\ol{G}(\e,\ol{\be}_{0}(\e))=0$. 
In particular $\be(t;\e)=\ol{\be}_{0}(\e)+b^{\RR}(t;\e,\ol{\be}_{0}(\e))$ 
defined in (\ref{eq:3.18}) 
solves the equation of motion (\ref{eq:1.1}). 
\end{lemma} 
 
\prova 
It follows from the results above. 
\EP 
 
\zerarcounters 
\section{Proof of Theorem \ref{thm:2.2}} 
\label{sec:6} 
 
If $F_{\vzero}(\be_{0})$ vanishes identically, let us come back to the 
formal expansion (\ref{eq:3.4}) of $G(\e,\be_{0})$, 
where $G^{(0)}(\be_{0})=F_{\vzero}(\be_{0})\equiv0$ by hypothesis. 
 
Assume first that there exists $k_{0}\in\NNN$ such that all functions 
${G}^{(k)}(\be_{0})$ are identically zero for $0\le k \le k_{0}-1$, while 
${G}^{(k_{0})}(\be_{0})$ is not identically vanishing. Then we can write 
\begin{equation}\label{eq:6.2a} 
{G}(\e,\be_{0})  = \e^{k_{0}} \left({G}^{(k_{0})} 
(\be_{0})+ {G}^{(> k_{0})}(\e,\be_{0}) \right), 
\end{equation} 
with $G^{(>k_{0})}(\e,\be_{0})=O(\e)$, and we can solve the equation 
of motion up to order $k_{0}$ without fixing the parameter $\be_{0}$. 
 
Any primitive function $g^{(k_{0})}(\be_{0})$ of ${G}^{(k_{0})}(\be_{0})$ 
is therefore analytic and periodic: since it is not identically constant, 
it admits at least one maximum $\bar{\be}_{0}'$ and one minimum 
$\bar{\be}_{0}''$, so that one can assume the following 
 
\begin{hyp}\label{hyp3} 
${\be}_{0}^{*}$ is a zero of order $\bar{\gotn}$ for 
${G}^{(k_{0})}(\be_{0})$ with $\bar{\gotn}$ odd, and 
$\e^{k_{0}+1} \partial^{\bar{\gotn}}_{\be_{0}}{G}^{(k_{0})} ({\be}_{0}^{*})<0$. 
\end{hyp} 
 
Indeed, if $k_{0}$ is even one can choose ${\be}^{*}_{0}=\bar{\be}_{0}'$ 
for $\e>0$, and ${\be}^{*}_{0}=\bar{\be}_{0}''$ for $\e<0$; if 
$k_{0}$ is odd we have to fix ${\be}^{*}_{0}=\bar{\be}_{0}'$: in both 
cases Hypothesis \ref{hyp3} is satisfied. 
 
Then one can adapt the proof in the previous sections to cover this case. 
Namely, as the formal expansion of $G^{\RR}$ coincide with that of $G$, 
one sets
$$ 
G^{\RR}(\e,\be_{0})=:\e^{k_{0}}G_{*}(\e,\be_{0}), 
$$ 
and hence, if $\MM$ satisfies property 1, 
\begin{equation}\label{eq:6.3} 
{\MM}_{\io}(0;\e,\be_{0})= 
\e^{k_{0}+1}\partial_{\be_{0}}G_{*}(\e,\be_{0}). 
\end{equation} 
On the other hand, Hypothesis \ref{hyp3} and Lemma 
\ref{lem:4.11} guarantee 
the existence of a continuous curve $\be_{0}(\e)$ such that
$\be_{0}(0)= \be_{0}^{*}$,
${G}_{*}(\e,\be_{0}(\e))\equiv 0$ and if $k_{0}$ is 
even then $\e^{k_{0}+1}\partial_{\be_{0}}{G}_{*}(\e,\be_{0}(\e)) 
\le0$ in a suitable half-neighbourhood of $\e=0$, while if $k_{0}$ is odd 
and ${\be}_{0}^{*}$ is a maximum for $g^{(k_{0})}$, then 
$\partial_{\be_{0}}{G}_{*}(\e,\be_{0}(\e))\le 0$ in a whole neighbourhood 
of $\e=0$. Then one can reason as in Section \ref{sec:5} to obtain the result. 
 
Finally, assume $G^{(k)}(\be_{0})\equiv 0$ for all $k\ge0$. 
We shall see that no resummation is necessary in that case: 
this situation is reminiscent of the ``null-renormalisation'' case 
considered in \cite{GCB} when studying the stability problem 
for Hill's equation with a quasi-periodic perturbation. 
 
We define trees and clusters according to the definitions previously done.
On the other hand, we slight change the definition of self-energy
clusters. Namely, a cluster $T$ on scale $n\ge0$ with only one entering 
line $\ell'_{T}$ and one exiting line $\ell_{T}$, and with $\nn_{\ell_{T}}= 
\nn_{\ell'_{T}}$, is called a self-energy cluster if $n+2\le n_{T}:=\min\{ 
n_{\ell_{T}}, n_{\ell'_{T}}\}$. The definition of self-energy cluster does not
change for the self-energy cluster on scale $-1$.
We denote by $\Theta_{k,\nn}$ the set of trees with 
order $k$ and momentum $\nn$ as in Section \ref{sec:3}, and by
$\gotS^{k}_{n}$ the set of 
(non-renormalised) self-energy clusters with order $k$ and scale $n$;
note that self-energy clusters
are allowed both in $\Theta_{k,\nn}$ and in $\gotS^{k}_{n}$.

For any subgraph $S$ of any tree $\theta\in\Theta_{k,\nn}$, and for any 
$T\in\gotS^{k}_{n}$, 
define the (non-renormalised) value of $S$ and $T$ 
as in (\ref{eq:3.15}) and (\ref{eq:3.14}) respectively, but with the 
(undressed) propagators defined as 
\begin{equation}\label{eq:6.2} 
\calG_{n_{\ell}}(\oo\cdot\nn_{\ell}):= 
\left\{\begin{aligned} 
&\frac{\Psi_{n_{\ell}}(\oo\cdot\nn_{\ell})}{\oo\cdot\nn_{\ell}^{2}}, 
   &n_{\ell}\ge 0,\\ 
&\; 1, &n_{\ell}=-1. 
\end{aligned}\right. 
\end{equation} 
Note that now the values of trees and self-energy clusters
do not depend on $\e$, and they depend on $\be_{0}$ only through the 
node factors. From now on we do not write explicitly the dependence 
on $\be_{0}$ to lighten the notations.  For all $k\ge1$, define
\begin{subequations} 
\begin{align} 
b^{(k)}_{\nn}& :=\sum_{\theta\in\Theta_{k,\nn}}\Val(\theta), 
\label{eq:6.3a}\\ 
G^{(k-1)}& :=\sum_{\theta\in\Theta_{k,\vzero}}\Val(\theta), 
\label{eq:6.3b}\\ 
M^{(k)}_{n}(x)& :=\sum_{T\in\gotS^{k}_{n}}\Val_{T}(x), 
\quad n\ge -1 
\label{eq:6.3c}\\ 
\MM^{(k)}_{n}(x)& :=\sum_{p=0}^{n}M^{(k)}_{p}(x), 
\quad n\ge -1 
\label{eq:6.3d}\\ 
\MM^{(k)}_{\io}(x)& :=\lim_{n\to\io} \MM^{(k)}_{n}(x). 
\label{eq:6.3e} 
\end{align} 
\label{eq:6.3bis} 
\end{subequations} 
The coefficients (\ref{eq:6.3a}) and (\ref{eq:6.3b}) 
coincide with (\ref{eq:3.2}) and (\ref{eq:3.5}) respectively as 
it is easy to check; in particular, for all $k\ge1$ one has 
$$ 
\sum_{\theta\in\Theta_{k,\vzero}}\Val(\theta) \equiv 0,
$$ 
by assumption. 
 
\begin{rmk}\label{rmk:6.1} 
\emph{ 
One has $\gotS^{k}_{-1}=\gotS^{1}_{n}=\emptyset$ for $k\ge2$ and $n\ge0$.
On the other hand $|\gotS^{1}_{-1}|=1$ and $\Val_{T}(x)=\partial_{\be_{0}}
F_{\vzero}\equiv0$ if $T$ is the self-energy cluster in $\gotS^{1}_{-1}$;
see Remark \ref{rmk:3.5bis}. Hence
$$ 
M^{(1)}_{n}(x)=\MM^{(1)}_{n}(x)=\MM^{(1)}_{\io}(x)
=M^{(k)}_{-1}=\MM^{(k)}_{-1}\equiv 0
$$ 
for all $n\ge-1$, $k\ge1$. 
} 
\end{rmk} 
 
Given a tree $\theta$ with $\Val(\theta)\ne 0$, we shall say that a line
$\ell\in L(\theta)$  is \emph{resonant} if it is the exiting line of a
self-energy  cluster $T$, otherwise we shall say that $\ell$ is
\emph{non-resonant}.
For any subgraph $T$ of any tree $\theta\in\Theta_{k,\nn}$, 
denote by $\gotN^{*}_{n}(T)$ the number 
of non-resonant lines on scale $\ge n$ in $T$, and set $K(T)$ 
as in (\ref{eq:3.22}). 
Then we can prove the analogous of Lemmas \ref{lem:3.4} and \ref{lem:3.5}, 
namely the following results. 
 
\begin{lemma}\label{lem:6.2} 
For any $\theta \in \Theta_{k,\nn}$ such that $\Val(\theta)\ne 0$ 
one has $\gotN^{*}_{n}(\theta)\le 2^{-(m_{n}-2)}K(\theta)$, for all $n\ge0$. 
\end{lemma} 
 
\begin{lemma}\label{lem:6.2b} 
For any $T \in \gotS^{k}_{n}$ such that $\Val_{T}(x)\ne 0$ 
one has $\gotN^{*}_{p}(T)\le 2^{-(m_{p}-2)}K(T)$, for all $0\le p\le n$. 
\end{lemma} 
 
We omit the proofs of the two results above as it would be essentially 
a repetition of those for Lemmas \ref{lem:3.4} and \ref{lem:3.5}, 
respectively. Note that, since self-energy clusters are now allowed, for 
the proof of Lemma \ref{lem:6.2b} one 
needs that the momenta of the lines in $\calP_{T}$ are different from those 
of the external lines: this explains the new definition of 
self-energy clusters. 
 
In light of Lemmas \ref{lem:6.2} and \ref{lem:6.2b}, although one has 
the `good 
bound' $1/x^{2}$ for the propagators, one cannot prove the convergence 
of the power series (\ref{eq:3.1}) as done in Lemma \ref{lem:4.2}, 
because we do not have any bound for the number of resonant lines, which
in principle can accumulate `too much'. In fact, we need a gain factor 
proportional to $(\oo\cdot\nn_{\ell})^{2}$ for each resonant line $\ell$. 
 
\begin{lemma}\label{lem:6.3} 
For all $n\ge0$ and for all $k\ge2$ one has 
$\partial_{x}M_{n}^{(k)}(0)=0$, and hence $\partial_{x}\MM_{n}^{(k)}(0)=0$
for all $k\ge2$. 
\end{lemma} 
 
\prova 
As the propagators are trivially even in the momenta, one can repeat 
(almost word by word) the proof of Lemma \ref{lem:4.5} so as to 
obtain the result. 
\EP 
 
\begin{lemma}\label{lem:6.4} 
One has $\MM_{\io}^{(k)}(0)\equiv 0$ for all $k\ge 2$. 
\end{lemma} 
 
\prova 
One has (see also Remark \ref{rmk:4.9}) $\partial_{\be_{0}}G^{(k-1)}\equiv 
\MM_{\io}^{(k)}(0)$ so that the assertion follows. 
\EP 
 
\begin{lemma}\label{lem:6.5} 
For all $k\ge 1$ one has 
\begin{equation}\label{eq:6.4} 
|\MM^{(k)}_{n}(x)|\Psi_{n+2}(x)\le C^{k}x^{2}\Psi_{n+2}(x), 
\end{equation} 
for some positive constant $C$. 
\end{lemma} 
 
 
\prova 
First of all note that (\ref{eq:6.4}) is trivially satisfied 
if $\Psi_{n+2}(x)=0$. Assume then 
\begin{equation}\label{eq:6.5} 
\frac{\al_{m_{n+2}}(\oo)}{8}<|x|<\frac{\al_{m_{n+1}}(\oo)}{4}. 
\end{equation} 
Note also that the bound (\ref{eq:6.4}) provides the gain factor which
is needed for the resonant lines. This can be seen as follows. 
 
Let $\theta \in \Theta_{k,\nn}$ and let $S$ be any subgraph of $\theta$. 
For any $\ell\in L(S)$ set 
\begin{equation}\label{eq:6.6bis} 
\calA_{\ell}(S,x_{\ell}):= 
\Big(\prod_{\substack{v\in N(S)\\ v \not\prec \ell}}\calF_{v}\Big) 
\Big(\prod_{\substack{\ell'\in L(S)\\ \ell' \not\preceq\ell}}\calG_{n_{\ell'}} 
(x_{\ell'})\Big), 
\end{equation} 
and 
\begin{equation}\label{eq:6.7bis} 
\BB_{\ell}(S):= 
\Big(\prod_{\substack{v\in N(S)\\ v\prec \ell}}\calF_{v}\Big) 
\Big(\prod_{\substack{\ell'\in L(S)\\ \ell' \prec\ell}}\calG_{n_{\ell'}} 
(x_{\ell'})\Big) , 
\end{equation} 
where $x_{\ell}=\oo\cdot\nn_{\ell}$. 
If $\ell\in L(S)$ is a resonant line exiting a self-energy cluster 
$T\in\gotS^{k}_{n}$ (with of course $n\le n_{T}-2\le n_{\ell}-2$) and also 
$\ell'_{T}\in L(S)$ we can write 
\begin{equation}\label{eq:6.8bis} 
\Val(S)=\calA_{\ell}(S,x_{\ell})\calG_{n_{\ell}}(x_{\ell}) 
 \Val_{T}(x_{\ell})\calG_{n_{\ell'_{T}}}(x_{\ell})\BB_{\ell'_{T}}(S), 
\end{equation} 
where we have used $x_{\ell}=x_{\ell'_{T}}$. 
But then, if we sum over all $S'$ which can be obtained from 
$S$ by replacing 
$T$ with any self-energy cluster $T'\in\gotS^{k}_{n'}$ for any 
$n'\le n_{T}-2$ we obtain 
\begin{equation}\label{eq:6.9bis} 
\calA_{\ell}(S,x_{\ell})\calG_{n_{\ell}}(x_{\ell}) 
 \MM_{n_{T}-2}^{(k)}(x_{\ell})\calG_{n_{\ell'_{T}}}(x_{\ell})\BB_{\ell'_{T}}(S), 
\end{equation} 
and hence, by (\ref{eq:6.4}), we obtain the gain factor which is needed. 
 
We shall prove the bound (\ref{eq:6.6}) by induction on $k$. 
For $k=1$ (\ref{eq:6.4}) is trivially satisfied. Assume (\ref{eq:6.4}) 
to hold for all $k'<k$. 
By Lemma \ref{lem:6.3} we can write 
\begin{equation}\label{eq:6.6} 
\MM_{n}^{(k)}(x)=\MM_{n}^{(k)}(0)+x^{2}\int_{0}^{1}dt\,(1-t)\partial^{2} 
\MM_{n}^{(k)}(tx), 
\end{equation} 
where $\partial^{2}$ denotes the second derivative of $\MM_{n}^{(k)}$ 
with respect to its argument. 
Then we shall prove 
\begin{subequations} 
\begin{align} 
|\MM^{(k)}_{n}(0)|&\le A_{1}^{k}\frac{\al_{m_{n+2}}^{2}(\oo)}{64}, 
\label{eq:6.7a}\\ 
|\partial^{2}\MM^{(k)}_{n}(x)|&\le A_{2}^{k}, 
\label{eq:6.7b} 
\end{align} 
\label{eq:6.7} 
\end{subequations} 
for suitable constants $A_{1},A_{2}$. 
Note that Lemma \ref{lem:6.2b} and the inductive 
hypothesis yield 
\begin{subequations} 
\begin{align} 
|M_{n}^{(k)}(x)|&\le B_{1}^{k}e^{-B_{2}2^{m_{n}}}, 
\label{eq:6.8a}\\ 
|\partial^{2}M_{n}^{(k)}(x)|&\le D_{1}^{k}e^{-D_{2}2^{m_{n}}}, 
\label{eq:6.8b} 
\end{align} 
\label{eq:6.8} 
\end{subequations} 
for some positive constants $B_{1},B_{2}, D_{1}$ and $D_{2}$. But then 
\begin{equation}\label{eq:6.9} 
|\MM_{n}^{(k)}(0)|=|\MM_{n}^{(k)}(0)-\MM_{\io}^{(k)}(0)|
\le\sum_{p\ge n+1}|M_{p}^{(k)}(0)| \le 2 B_{1}^{k}e^{-B_{2}2^{m_{n}}}, 
\end{equation} 
so that (\ref{eq:6.7a}) follows if $A_{1}$ is suitably chosen. 
Moreover 
\begin{equation}\label{eq:6.11} 
|\partial^{2}\MM_{n}^{(k)}(x)|\le\sum_{p=0}^{n}|\partial^{2}M_{p}^{(k)}(x)|\le
D D_{1}^{k} 
\end{equation} 
for some constant $D$. Hence the assertion follows. 
\EP 
 
\begin{rmk}\label{rmk:6.7} 
\emph{
We have obtained the convergence of the power series (\ref{eq:3.1}) 
and (\ref{eq:3.4}) for any $\be_{0}$ and any $\e$ small enough. 
Hence, in this case, the response solution turns out to be analytic in 
both $\e,\be_{0}$. 
}
\end{rmk} 

\begin{rmk}\label{rmk:6.8} 
\emph{
Note that the problem under study has analogies with the problem 
considered in \cite{G4}. In that case, the resummation adds to the 
small divisor $\ii\oo\cdot\nn$ a quantity 
$-\e(\oo\cdot\nn)^{2}+\MM_{n}(\oo\cdot\nn;\e)$, and one can prove that 
$\MM_{n}(x,\e)$ is smooth in $x$ and it is real at $x=0$, 
so that the dressed propagator is proportional to 
$1/(\ii\oo\cdot\nn-\e(\oo\cdot\nn)^{2}+\MM_{n}(\oo\cdot\nn;\e))$, 
and hence can be bounded essentially as the undressed one. 
In the present case, both the small divisor $(\oo\cdot\nn)^{2}$ 
and the correction are real, but they turn 
out to have the same sign (for a suitable choice of $\be_{0}^{*}$), 
so that once more the dressed propagator can be bounded 
as the undressed one. 
}
\end{rmk} 
 
\appendix 
 
\zerarcounters 
\section{Proof of Lemma \ref{lem:4.8}} 
\label{app:a} 
 
First of all, for any renormalised tree $\theta$ set 
\begin{equation} \label{eq:a.1} 
\partial_{v}\Val(\theta;\e,\be_{0}) := 
\partial_{\be_{0}}\calF_{v}(\be_{0})\left(\prod_{w \in N(\theta)\setminus\{v\}} 
\calF_{w}(\be_{0})\right) 
\left(\prod_{\ell\in L(\theta)}\calG_{n_{\ell}}(\oo\cdot\nn_{\ell}; 
\e,\be_{0})\right) 
\end{equation} 
and 
\begin{equation}\label{eq:a.2} 
\begin{aligned} 
\partial_{\ell}\Val(\theta;\e,\be_{0})&:= 
\partial_{\be_{0}}\calG_{n_{\ell}}(x_{\ell};\e,\be_{0}) 
\left(\prod_{v\in N(\theta)}\calF_{v}(\be_{0}) 
\right) \left(\prod_{\la\in L(\theta)\setminus\{\ell\}} 
\calG_{n_{\la}}(x_{\la};\e,\be_{0})\right)\\ 
 &=\calA_{\ell} (\theta,x_{\ell};\e,\be_{0}) \, 
\partial_{\be_{0}}\calG_{n_{\ell}}(x_{\ell};\e,\be_{0})\, 
\BB_{\ell}(\theta;\e,\be_{0}), 
\end{aligned} 
\end{equation} 
where $x_{\ell}:=\oo\cdot\nn_{\ell}$, 
$\partial_{\be_{0}}\calG_{n_{\ell}}(x_{\ell};\e,\be_{0})$ 
is written according to Remark \ref{rmk:3.1ter}, 
\begin{equation}\label{eq:a.3} 
\calA_{\ell}(\theta,x_{\ell};\e,\be_{0}):= 
\Big(\prod_{\substack{v\in N(\theta)\\ v \not\prec \ell}}\calF_{v}(\be_{0})\Big) 
\Big(\prod_{\substack{\ell'\in L(\theta)\\ \ell' \not\preceq\ell}}\calG_{n_{\ell'}} 
(x_{\ell'};\e,\be_{0})\Big), 
\end{equation} 
and 
\begin{equation}\label{eq:a.4} 
\BB_{\ell}(\theta;\e,\be_{0}):= 
\Big(\prod_{\substack{v\in N(\theta)\\ v\prec \ell}}\calF_{v}(\be_{0})\Big) 
\Big(\prod_{\substack{\ell'\in L(\theta)\\ \ell' \prec\ell}}\calG_{n_{\ell'}} 
(x_{\ell'};\e,\be_{0})\Big) , 
\end{equation} 
see also (\ref{eq:6.6bis}) and (\ref{eq:6.7bis}). 
Let us define in the analogous way $\partial_{v}\Val_{T}(x;\e,\be_{0})$ 
and $\partial_{\ell}\Val_{T}(x;\e,\be_{0})$ for any self-energy cluster $T$, 
and let us write 
\begin{equation} \label{eq:a.5} 
\partial_{\be_{0}}\Val(\theta;\e,\be_{0})= 
\partial_{N}\Val(\theta;\e,\be_{0}) + 
\partial_{L}\Val(\theta;\e,\be_{0}), 
\end{equation} 
where 
\begin{equation} \label{eq:a.6} 
\partial_{N}\Val(\theta;\e,\be_{0}) :=\sum_{v \in N(\theta)} 
\partial_{v}\Val(\theta;\e,\be_{0}), 
\end{equation} 
and 
\begin{equation}\label{eq:a.7} 
\begin{aligned} 
\partial_{L}\Val(\theta;\e,\be_{0})&:=\sum_{\ell \in L(\theta)} 
\partial_{\ell}\Val(\theta;\e,\be_{0}). \\ 
\end{aligned} 
\end{equation} 
Let us also write 
\begin{equation}\label{eq:a.8} 
\partial_{\be_{0}}\Val_{T}(x;\e,\be_{0})= 
\partial_{N}\Val_{T}(x;\e,\be_{0})+\partial_{L}\Val_{T}(x;\e,\be_{0}), 
\end{equation} 
for any $T\in\gotR_{n}$, $n\ge0$, where the derivatives $\partial_{N}$ 
and $\partial_{L}$ are defined analogously with the previous cases 
(\ref{eq:a.6}) and (\ref{eq:a.7}), with $N(T)$ and $L(T)$ replacing 
$N(\theta)$ and $L(\theta)$, respectively, so that we can split 
\begin{equation}\label{eq:a.9} 
\begin{aligned} 
\partial_{\be_{0}}M_{n}(x;\e,\be_{0})&= 
\partial_{N}M_{n}(x;\e,\be_{0})+\partial_{L}M_{n}(x;\e,\be_{0}), \\ 
\partial_{\be_{0}}\MM_{n}(x;\e,\be_{0})&= 
\partial_{N}\MM_{n}(x;\e,\be_{0})+\partial_{L}\MM_{n}(x;\e,\be_{0}) , 
\end{aligned} 
\end{equation} 
again with obvious meaning of the symbols. 
 
\begin{rmk}\label{rmk:a.1} 
\emph{ 
We can interpret the derivative $\partial_{v}$ as all the possible 
ways to attach an extra line (carrying a momentum $\vzero$) to 
the node $v$, so that 
$$ 
\sum_{k\ge0}\e^{k+1}\sum_{\theta\in\Theta^{\RR}_{k+1,\vzero}}\partial_{N} 
\Val(\theta;\e,\be_{0}), 
$$ 
produces contributions to $\MM_{\io}(0;\e,\be_{0})$. 
} 
\end{rmk} 
 
Given any $\theta \in \Theta^{\RR}_{k,\vzero}$ we have to study the 
derivative (\ref{eq:a.5}). The terms (\ref{eq:a.6}) produce immediately 
contributions to $\MM_{\io}(0;\e,\be_{0})$ by Remark \ref{rmk:a.1}. 
Thus, we have to study the derivatives $\partial_{\ell}\Val(\theta;\e,\be_{0})$ 
appearing in the sum (\ref{eq:a.7}). 
Here and henceforth, we shall not write any longer explicitly the dependence 
on $\e$ and $\be_{0}$ of both propagators and self-energies, 
in order not to overwhelm the notation. 
 
For any $\theta \in \Theta^{\RR}_{k,\vzero}$ 
such that $\Val(\theta;\e,\be_{0})\ne 0$ and for any line $\ell\in L(\theta)$, 
either there is only one scale $n$ such that $\Psi_{n}(x_{\ell}) \neq 0$ 
(and in that case $\Psi_{n}(x_{\ell})=1$ and $\Psi_{n'}(x_{\ell})=0$ 
for all $n'\neq n$) or there exists only one 
$n\ge0$ such that $\Psi_{n}(x_{\ell})\Psi_{n+1}(x_{\ell})\ne 0$. 
 
\noindent 
\textbf{1.} If $\Psi_{n}(x_{\ell})=1$ one has 
\begin{equation}\label{eq:a.10} 
\begin{aligned} 
\partial_{\ell}\Val(\theta;\e,\be_{0})&=\calA_{\ell}(\theta,x_{\ell}) \, 
\frac{\Psi_{n}(x_{\ell})}{x_{\ell}^{2}- 
\MM_{n-1}(x_{\ell})} 
\partial_{\be_{0}}\MM_{n-1}(x_{\ell}) 
\frac{1}{x_{\ell}^{2}-\MM_{n-1}(x_{\ell})} \, \BB_{\ell}(\theta)\\ 
&=\calA_{\ell}(\theta,x_{\ell}) \, 
\frac{\Psi_{n}(x_{\ell})}{x_{\ell}^{2}- 
\MM_{n-1}(x_{\ell})} 
\partial_{\be_{0}}\MM_{n-1}(x_{\ell}) 
\frac{\Psi_{n}(x_{\ell})}{x_{\ell}^{2}-\MM_{n-1}(x_{\ell})} 
\, \BB_{\ell}(\theta) \\ 
&=\calA_{\ell}(\theta,x_{\ell}) \, 
\calG_{n}(x_{\ell}) 
\partial_{\be_{0}}\MM_{n-1}(x_{\ell}) 
\calG_{n}(x_{\ell}) \, \BB_{\ell}(\theta) , 
\end{aligned} 
\end{equation} 
where (here and henceforth) we shorten $\calA_{\ell}(\theta,x_{\ell}) 
=\calA_{\ell}(\theta,x_{\ell};\e,\be_{0})$ and $\BB_{\ell}(\theta) 
=\BB_{\ell}(\theta;\e,\be_{0})$. 
 
\begin{rmk}\label{rmk:a.2} 
\emph{ 
Note that if we split $\partial_{\be_{0}}=\partial_{N}+ 
\partial_{L}$ in (\ref{eq:a.10}), the term with $\partial_{N}\MM_{n-1} 
(x_{\ell})$ is a contribution to $\MM_{\io}(0)$. 
} 
\end{rmk} 
 
If there is only one $n\ge0$ such that $\Psi_{n}(x_{\ell})\Psi_{n+1}(x_{\ell}) 
\ne 0$, then $\Psi_{n}(x_{\ell})+\Psi_{n+1}(x_{\ell})=1$ and $\chi_{q}(x_{\ell})=1$ for all 
$q=-1,\ldots,n-1$, so that $\psi_{n+1}(x_{\ell})=1$ and hence 
$\Psi_{n+1}(x_{\ell})=\chi_{n}(x_{\ell})$. Moreover
it can happen only (see Remark \ref{rmk:3.2}) $n_{\ell}=n$ or
$n_{\ell}=n+1$. 

\begin{figure}[ht] 
\centering 
\ins{040pt}{-40pt}{$n\!+\!1$} 
\ins{096pt}{-47pt}{$\le n$} 
\ins{144pt}{-40pt}{$n\!+\!1$} 
\ins{224pt}{-48pt}{$+$} 
\ins{288pt}{-40pt}{$n\!+\!1$} 
\ins{332pt}{-47pt}{$\le n\!-\!1$} 
\ins{406pt}{-41pt}{$n$} 
\ins{080pt}{-113pt}{$+$} 
\ins{142pt}{-105pt}{$n\!+\!1$} 
\ins{190pt}{-112pt}{$\le n\!-\!1$} 
\ins{250pt}{-106pt}{$n$} 
\ins{308pt}{-113pt}{$n$} 
\ins{353pt}{-105pt}{$n\!+\!1$} 
\includegraphics[width=6in]{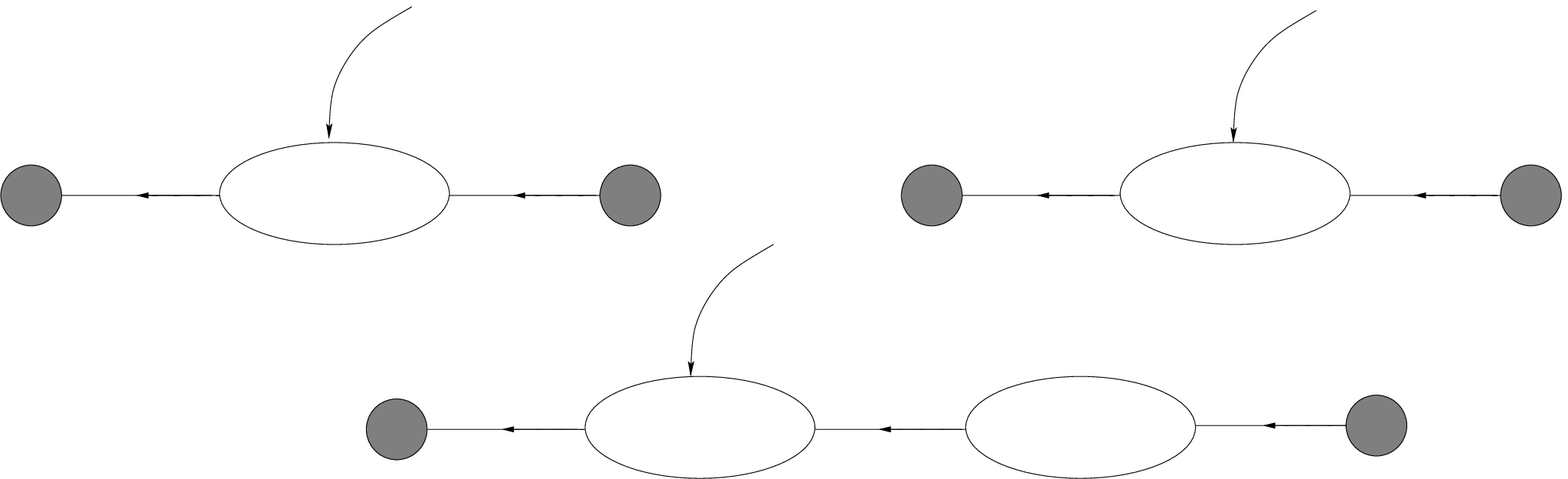} 
\caption{Graphical representation of the derivative 
$\partial_{\ell}\Val(\theta;\e,\be_{0})$ according to (\ref{eq:a.11}).} 
\label{fig:5} 
\end{figure}

\begin{figure}[ht] 
\centering 
\ins{047pt}{-41pt}{$n$} 
\ins{088pt}{-47pt}{$\le n-1$} 
\ins{152pt}{-41pt}{$n$} 
\ins{224pt}{-48pt}{$+$} 
\ins{300pt}{-41pt}{$n$} 
\ins{332pt}{-47pt}{$\le n\!-\!1$} 
\ins{390pt}{-40pt}{$n\!+\!1$} 
\ins{080pt}{-113pt}{$+$} 
\ins{150pt}{-106pt}{$n$} 
\ins{190pt}{-112pt}{$\le n\!-\!1$} 
\ins{255pt}{-106pt}{$n$} 
\ins{308pt}{-113pt}{$n$} 
\ins{353pt}{-105pt}{$n\!+\!1$} 
\includegraphics[width=6in]{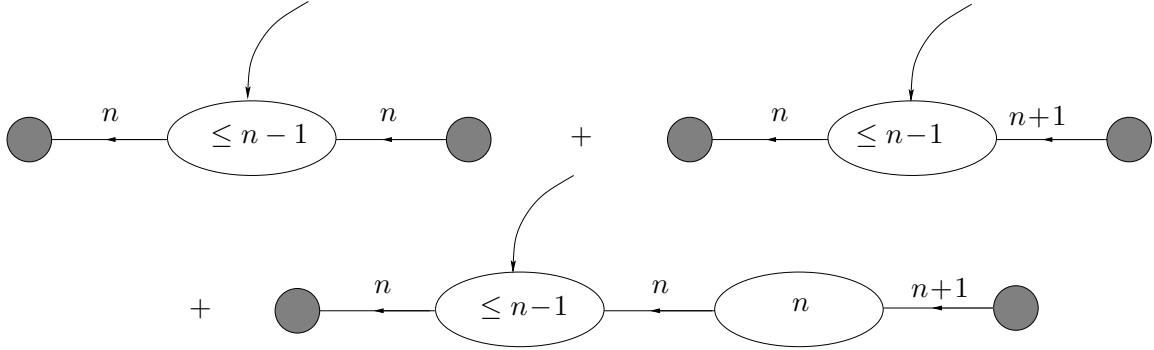} 
\caption{Graphical representation of the three contributions in 
the last two lines of (\ref{eq:a.13}).} 
\label{fig:7} 
\end{figure}

\begin{figure}[ht] 
\centering 
\ins{030pt}{-28pt}{$n\!+\!1$} 
\ins{063pt}{-38pt}{$n$} 
\ins{095pt}{-29pt}{$n$} 
\ins{120pt}{-36pt}{$\le \! n\!-\!1$} 
\ins{175pt}{-29pt}{$n$} 
\ins{222pt}{-38pt}{$+$} 
\ins{262pt}{-28pt}{$n\!+\!1$} 
\ins{299pt}{-38pt}{$n$} 
\ins{328pt}{-29pt}{$n$} 
\ins{356pt}{-36pt}{$\le n\!-\!1$} 
\ins{408pt}{-28pt}{$n\!+\!1$} 
\ins{066pt}{-088pt}{$+$} 
\ins{109pt}{-080pt}{$n\!+\!1$} 
\ins{143pt}{-090pt}{$n$} 
\ins{175pt}{-081pt}{$n$} 
\ins{201pt}{-087pt}{$\le n\!-\!1$} 
\ins{256pt}{-081pt}{$n$} 
\ins{297pt}{-089pt}{$n$} 
\ins{330pt}{-080pt}{$n\!+\!1$} 
\includegraphics[width=6.0in]{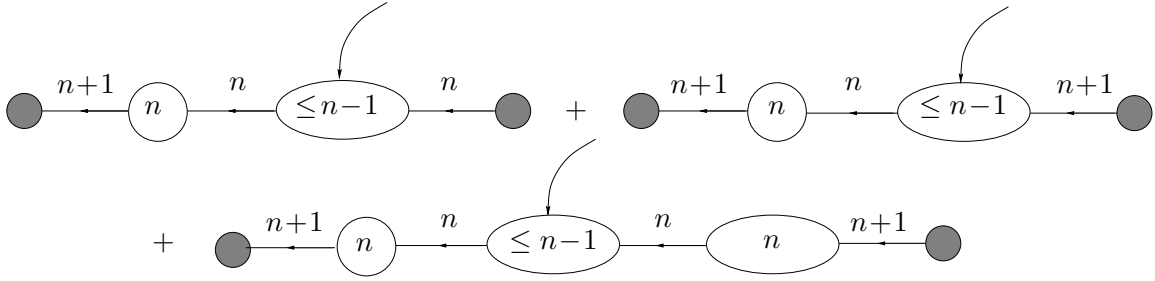} 
\caption{Graphical representation of the three contributions arising from 
(\ref{eq:a.16}).} 
\label{fig:10} 
\end{figure}

\begin{figure}[ht] 
\vskip-.3truecm 
\centering 
\ins{055pt}{-49pt}{$n\!+\!i$} 
\ins{060pt}{-64pt}{$\ell_{i}$} 
\ins{090pt}{-56pt}{$\le\!\! n\!\!+\!\!i\!\!-\!\!1$} 
\ins{118pt}{-80pt}{$S_{i-1}$} 
\ins{130pt}{-26pt}{$R_{i}$} 
\ins{135pt}{-64pt}{$\ell_{i-1}$} 
\ins{170pt}{-56pt}{$\le\!\! n\!\!+\!\!i\!\!-\!\!1$} 
\ins{160pt}{-80pt}{$S'_{i-1}$} 
\ins{226pt}{-64pt}{$\ell'_{i}$} 
\ins{220pt}{-49pt}{$n\!+\!i$} 
\ins{256pt}{-54pt}{$n\!+\!i$} 
\ins{290pt}{-49pt}{$n\!+\!i\!+\!1$} 
\ins{300pt}{-64pt}{$\ell'_{_{T_{i}}}$} 
\ins{374pt}{-56pt}{$+$} 
\ins{055pt}{-153pt}{$n\!+\!i$} 
\ins{060pt}{-167pt}{$\ell_{i}$} 
\ins{091pt}{-158pt}{$\le\!\! n\!\!+\!\!i\!\!-\!\!1$} 
\ins{118pt}{-184pt}{$S_{i-1}$} 
\ins{135pt}{-167pt}{$\ell_{i-1}$} 
\ins{170pt}{-161pt}{$n\!+\!i\!-\!1$} 
\ins{220pt}{-153pt}{$n\!+\!i$} 
\ins{220pt}{-167pt}{$\ell'_{T_{i-1}}$} 
\ins{250pt}{-158pt}{$\le\!\! n\!\!+\!\!i\!\!-\!\!1$} 
\ins{277pt}{-184pt}{$S'_{i-1}$} 
\ins{295pt}{-153pt}{$n\!+\!i$} 
\ins{300pt}{-167pt}{$\ell'_{i}$} 
\ins{334pt}{-158pt}{$n+i$} 
\ins{372pt}{-153pt}{$n\!+\!i\!+\!1$} 
\ins{380pt}{-167pt}{$\ell'_{_{T_{i}}}$} 
\includegraphics[width=5.4in]{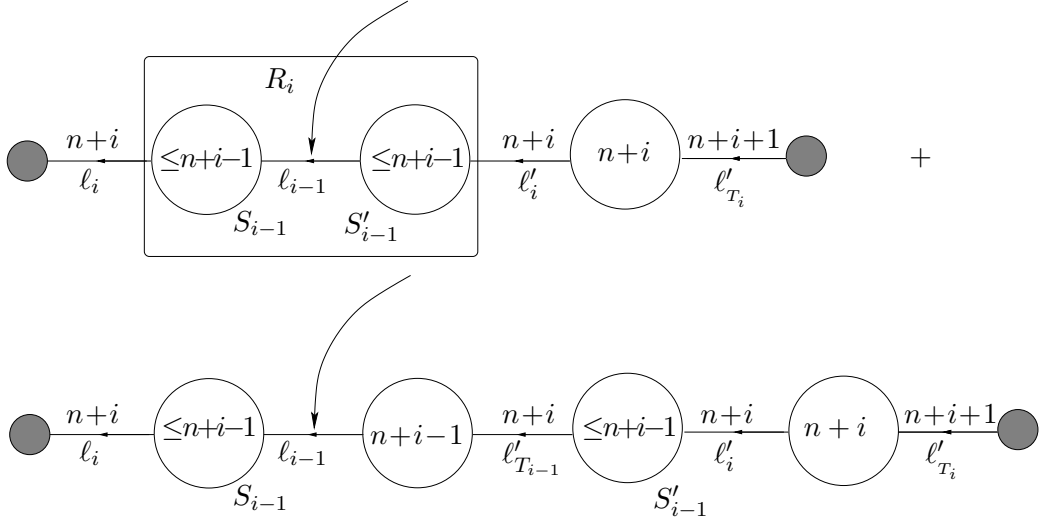} 
\vskip.3truecm 
\caption{Graphical representation of the left hand side of (\ref{eq:a.23}).} 
\label{fig:12} 
\end{figure}

\noindent 
\textbf{2.} Consider first the case $n_{\ell}=n+1$. One has
\begin{equation}\label{eq:a.11} 
\begin{aligned} 
\partial_{\ell}\Val(\theta;\e,\be_{0})& = 
  \calA_{\ell}(\theta,x_{\ell}) \, 
  \calG_{n+1}(x_{\ell}) 
  \partial_{\be_{0}}\MM_{n}(x_{\ell}) 
  \frac{1}{x_{\ell}^{2}- \MM_{n}(x_{\ell})} 
  \, \BB_{\ell}(\theta)\\ 
&=\calA_{\ell}(\theta,x_{\ell}) \, 
  \calG_{n+1}(x_{\ell}) 
  \partial_{\be_{0}}\MM_{n-1}(x_{\ell}) 
  \frac{\Psi_{n}(x_{\ell})+\Psi_{n+1}(x_{\ell})}{x_{\ell}^{2}- 
    \MM_{n}(x_{\ell})} 
  \, \BB_{\ell}(\theta)\\ 
&\qquad+\calA_{\ell}(\theta,x_{\ell}) \, 
  \calG_{n+1}(x_{\ell}) 
  \partial_{\be_{0}}M_{n}(x_{\ell}) 
  \frac{\chi_{n}(x_{\ell})}{x_{\ell}^{2}-\MM_{n}(x_{\ell})} 
  \, \BB_{\ell}(\theta) \\ 
&=\calA_{\ell}(\theta,x_{\ell}) \,  
  \calG_{n+1}(x_{\ell}) 
  \left(\sum_{q=-1}^{n}\partial_{\be_{0}}M_{q}(x_{\ell})\right) 
  \calG_{n+1}(x_{\ell}) 
  \, \BB_{\ell}(\theta)\\ 
&\qquad+\calA_{\ell}(\theta,x_{\ell}) \, 
  \calG_{n+1}(x_{\ell}) 
  \left(\sum_{q=-1}^{n-1}\partial_{\be_{0}}M_{q}(x_{\ell})\right) 
  \calG_{n}(x_{\ell}) 
  \, \BB_{\ell}(\theta)\\ 
&\qquad+\calA_{\ell}(\theta,x_{\ell}) \, 
  \calG_{n+1}(x_{\ell}) \!\! 
  \left(\sum_{q=-1}^{n-1}\partial_{\be_{0}}M_{q}(x_{\ell})\right) \!\! 
  \calG_{n}(x_{\ell}) M_{n}(x_{\ell}) 
  \calG_{n+1}(x_{\ell}) 
  \, \BB_{\ell}(\theta) . 
\end{aligned} 
\end{equation} 
We can represent graphically the three contributions in 
(\ref{eq:a.11}) as in Figure \ref{fig:5}: we represent the derivative
$\partial_{\be_{0}}$ as an arrow pointing toward the graphical representation
of the differentiated quantity; see also Figures \ref{fig:7}, \ref{fig:10}
and \ref{fig:12}. 
 
\begin{rmk}\label{rmk:a.3} 
\emph{ 
Note that the $M_{n}(x_{\ell})$ appearing in the latter line of 
(\ref{eq:a.11}) has 
to be interpreted (see Remark \ref{rmk:3.2}) as 
$$ 
\sum_{T\in\gotL\gotF_{n}}\e^{k(T)}\Val_{T}(x_{\ell};\e,\be_{0}). 
$$ 
Note also that, again, if we split $\partial_{\be_{0}}=\partial_{N}+ 
\partial_{L}$ in (\ref{eq:a.11}), all the terms with $\partial_{N}M_{q} 
(x_{\ell})$ are contributions to $\MM_{\io}(0)$. 
} 
\end{rmk} 
 
Now consider the case $n_{\ell}=n$. 
 
\begin{figure}[ht] 
\centering 
\ins{018pt}{-11pt}{$\bar{\theta}=$} 
\ins{076pt}{-04pt}{$n$} 
\ins{076pt}{-20pt}{$\ell$} 
\ins{193pt}{-11pt}{$\theta'=$} 
\ins{260pt}{-04pt}{$n$} 
\ins{260pt}{-20pt}{$\ell$} 
\ins{312pt}{-13pt}{$n$} 
\ins{368pt}{-03pt}{$n\!+\!1$} 
\includegraphics[width=5.2in]{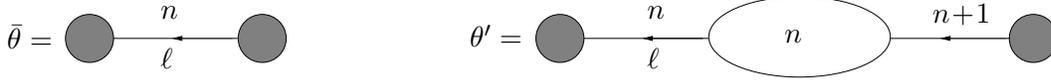} 
\caption{The renormalised tree $\bar{\theta}$ and the renormalised trees 
$\theta'$ of the set 
$\tau_{1}(\bar{\theta},\ell)$ associated with $\bar{\theta}$.} 
\label{fig:6} 
\end{figure}

\noindent 
\textbf{3.} If $\ell$ is not the exiting line of a 
left-fake cluster, set $\bar{\theta}=\theta$; otherwise, 
if $\ell$ is the exiting line of a left-fake cluster $T$, 
define -- if possible -- $\bar{\theta}$ as the renormalised tree 
obtained from $\theta$ by removing $T$ and $\ell_{T}'$.
In both cases, define -- if possible -- $\tau_{1}(\bar{\theta}, 
\ell)$ as the set constituted by all the renormalised trees $\theta'$ 
obtained from $\bar{\theta}$ by inserting a left-fake cluster, together 
with its entering line, between $\ell$ and the node $v$ which $\ell$ exits; 
see Figure \ref{fig:6}.
Here and henceforth, if $S$ is a subgraph with only one entering line
$\ell'_{S}=\ell_{v}$ and one exiting line $\ell_{S}$ and we ``remove'' $S$
together with $\ell'_{S}$, we mean that we also reattach the line $\ell_{S}$
to the node $v$.
 
\begin{rmk}\label{rmk:a.4} 
\emph{ 
The construction of the set $\tau_{1}(\bar{\theta},\ell)$ could be impossible 
if the removal or the insertion of a left-fake cluster $T$, together with its 
entering line $\ell'_{T}$, produce a self-energy cluster. We shall see 
later how to deal with these cases. 
} 
\end{rmk} 
 
Then one has 
\begin{equation}\label{eq:a.12} 
\partial_{\ell}\Val(\bar{\theta};\e,\be_{0})+ 
\partial_{\ell}\!\!\!\sum_{\theta'\in\tau_{1}(\bar{\theta},\ell)}\!\!\! 
 \Val(\theta';\e,\be_{0}) 
=\calA_{\ell}(\bar{\theta},x_{\ell}) \, \partial_{\be_{0}}\calG_{n}(x_{\ell}) 
\left(1+M_{n}(x_{\ell})\calG_{n+1}(x_{\ell})\right) \, \BB_{\ell}(\bar{\theta}), 
\end{equation} 
where 
\begin{equation}\label{eq:a.13} 
\begin{aligned} 
\partial_{\be_{0}}&\calG_{n}(x_{\ell})\left(1+ 
  M_{n}(x_{\ell})\calG_{n+1}(x_{\ell})\right)\\ 
&=\calG_{n}(x_{\ell}) 
  \partial_{\be_{0}}\MM_{n-1}(x_{\ell}) 
  \calG_{n}(x_{\ell})\\ 
&\quad+\calG_{n}(x_{\ell}) 
  \partial_{\be_{0}}\MM_{n-1}(x_{\ell}) 
  \frac{\Psi_{n+1}(x_{\ell})}{x_{\ell}^{2}-\MM_{n-1}(x_{\ell})}\\ 
&\quad+\calG_{n}(x_{\ell}) 
  \partial_{\be_{0}}\MM_{n-1}(x_{\ell}) 
  \calG_{n}(x_{\ell}) 
  M_{n}(x_{\ell}) 
  \calG_{n+1}(x_{\ell})\\  
&\quad+\calG_{n}(x_{\ell}) 
  \partial_{\be_{0}}\MM_{n-1}(x_{\ell}) 
  \frac{\Psi_{n+1}(x_{\ell})}{x_{\ell}^{2}-\MM_{n-1}(x_{\ell})} 
  M_{n}(x_{\ell}) 
  \calG_{n+1}(x_{\ell})\\  
&=\calG_{n}(x_{\ell}) 
  \partial_{\be_{0}}\MM_{n-1}(x_{\ell}) 
  \calG_{n}(x_{\ell}) +\calG_{n}(x_{\ell}) 
  \partial_{\be_{0}}\MM_{n-1}(x_{\ell}) 
  \calG_{n+1}(x_{\ell})\\ 
&\quad-\calG_{n}(x_{\ell}) 
  \partial_{\be_{0}}\MM_{n-1}(x_{\ell}) 
  \frac{\chi_{n}(x_{\ell})}{x_{\ell}^{2}-\MM_{n-1}(x_{\ell})} 
  M_{n}(x_{\ell}) 
  \calG_{n+1}(x_{\ell})\\  
&\quad+\calG_{n}(x_{\ell}) 
  \partial_{\be_{0}}\MM_{n-1}(x_{\ell}) 
  \calG_{n}(x_{\ell}) 
  M_{n}(x_{\ell}) 
  \calG_{n+1}(x_{\ell})\\ 
&\quad+\calG_{n}(x_{\ell}) 
  \partial_{\be_{0}}\MM_{n-1}(x_{\ell}) 
  \frac{\Psi_{n+1}(x_{\ell})}{x_{\ell}^{2}-\MM_{n-1}(x_{\ell})} 
  M_{n}(x_{\ell}) 
  \calG_{n+1}(x_{\ell})\\ 
&=\calG_{n}(x_{\ell}) 
  \partial_{\be_{0}}\MM_{n-1}(x_{\ell}) 
  \calG_{n}(x_{\ell}) +\calG_{n}(x_{\ell}) 
  \partial_{\be_{0}}\MM_{n-1}(x_{\ell}) 
  \calG_{n+1}(x_{\ell})\\ 
&\quad+\calG_{n}(x_{\ell}) 
  \partial_{\be_{0}}\MM_{n-1}(x_{\ell}) 
  \calG_{n}(x_{\ell}) 
  M_{n}(x_{\ell}) 
  \calG_{n+1}(x_{\ell}), 
\end{aligned} 
\end{equation} 
so that also in this case, if we split $\partial_{\be_{0}}=\partial_{N}+ 
\partial_{L}$, all the terms with $\partial_{N}\MM_{n-1}$ are 
contributions to $\MM_{\io}(0)$ -- see Remark \ref{rmk:a.2}. 
Again, we can represent graphically the three contributions obtained 
inserting (\ref{eq:a.13}) in (\ref{eq:a.12}): see Figure \ref{fig:7}. 
 
\begin{figure}[ht] 
\centering 
\ins{018pt}{-11pt}{$\theta'=$} 
\ins{072pt}{-03pt}{$n\!+\!1$} 
\ins{080pt}{-19pt}{$\ol{\ell}$} 
\ins{134pt}{-12pt}{$n$} 
\ins{192pt}{-04pt}{$n$} 
\ins{192pt}{-20pt}{$\ell$} 
\ins{310pt}{-11pt}{$\ol{\theta}=$} 
\ins{374pt}{-20pt}{$\ol{\ell}$} 
\ins{368pt}{-03pt}{$n\!+\!1$} 
\includegraphics[width=5.2in]{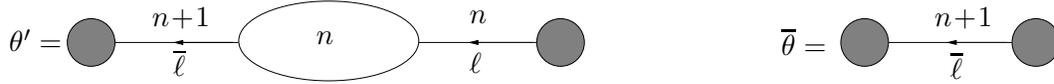} 
\caption{The  trees $\theta'$ of the set $\tau_{2}(\ol{\theta},\ol\ell)$ 
obtained from $\ol{\theta}$ when $\ell\in L(\theta)$ enters a right-fake 
cluster.} 
\label{fig:8} 
\end{figure}

\noindent 
\textbf{4.} Assume now that $\ell$ is not the exiting line of a left-fake cluster, 
and the insertion of a left-fake cluster, together 
with its entering line, produces a self-energy cluster. Note that this can 
happen only if $\ell$ is the entering line 
of a renormalised right-fake cluster $T$. Let $\ol{\ell}$ be 
the exiting line (on scale $n+1$) of the renormalised 
right-fake cluster $T$, call $\ol{\theta}$ the renormalised tree obtained 
from $\theta$ by removing $T$ and $\ell$ and call 
$\tau_{2}(\ol{\theta},\ol\ell)$ the set of renormalised trees $\theta'$ 
obtained from $\ol{\theta}$ by inserting a right-fake cluster, 
together with its entering line, before $\ol{\ell}$; see Figure \ref{fig:8}. 
 
By construction one has 
\begin{equation} \nonumber 
\begin{aligned} 
\Val(\ol{\theta};\e,\be_{0}) & = 
 \calA_{\ol{\ell}}(\ol{\theta},x_{\ell}) \, 
  \calG_{n+1}(x_{\ol{\ell}}) 
  \, \BB_{\ol{\ell}}(\ol{\theta}) \\ 
\sum_{\theta'\in\tau_{2}(\ol{\theta},\ol\ell)} 
\!\!\!\Val(\theta';\e,\be_{0}) & = 
 \calA_{\ol{\ell}}(\ol{\theta},x_{\ell}) \, 
  \calG_{n+1}(x_{\ol{\ell}}) \, M_{n}(x_{\ol{\ell}})\, 
  \calG_{n}(x_{\ol{\ell}}) 
  \, \BB_{\ol{\ell}}(\ol{\theta}), 
\end{aligned} 
\end{equation} 
where we have used that $x_{\ell}=x_{\bar{\ell}}$. 
 
\begin{figure}[ht] 
\centering 
\ins{183pt}{-20pt}{$T$} 
\ins{261pt}{-20pt}{$T''$} 
\ins{320pt}{-80pt}{$T'$} 
\ins{094pt}{-38pt}{$n\!+\!1$} 
\ins{168pt}{-45pt}{$n$} 
\ins{224pt}{-39pt}{$n$} 
\ins{282pt}{-45pt}{$n$} 
\ins{348pt}{-38pt}{$n\!+\!1$} 
\ins{226pt}{-54pt}{$\ell$} 
\ins{100pt}{-54pt}{$\ell_{T'}$} 
\ins{350pt}{-54pt}{$\ell_{T'}'$} 
\includegraphics[width=4.0in]{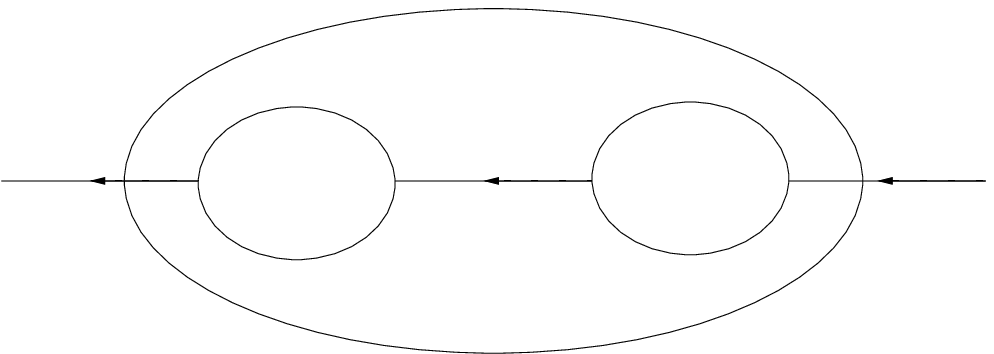} 
\caption{A self-energy cluster $T'\in \gotR_{n}(T)$.} 
\label{fig:9} 
\end{figure}

Consider the contribution to 
$\partial_{\ol{\ell}}\Val(\ol{\theta};\e,\be_{0})$ 
-- see (\ref{eq:a.11}) -- given by 
\begin{equation} \label{eq:a.14} 
\calA_{\ol{\ell}}(\ol{\theta},x_{\ol{\ell}}) 
  \calG_{n+1}(x_{\ol{\ell}}) 
  \partial_{L}M_{n}(x_{\ol{\ell}}) 
  \calG_{n+1}(x_{\ol{\ell}}) 
  \BB_{\ol{\ell}}(\ol{\theta}). 
\end{equation} 
Call $\gotR_{n}(T)$ the subset of $\gotR_{n}$ such that 
if $T'\in\gotR_{n}(T)$ the exiting line $\ell_{T'}$ exits also 
the renormalised right-fake cluster $T$; note that the entering line 
$\ell$ of $T$ must be also the exiting line of some renormalised 
left-fake cluster $T''$ contained in $T'$; see Figure \ref{fig:9}. 
 
Define 
\begin{equation} \label{eq:a.15} 
M_{n}(T,x_{\ol{\ell}};\e,\be_{o})= 
\sum_{T'\in\gotR_{n}(T)}\Val_{T'}(x_{\ol{\ell}};\e,\be_{0}). 
\end{equation} 
Hence one has 
\begin{equation}\label{eq:a.16} 
\begin{aligned} 
& \partial_{\ell} \!\!\!\sum_{\theta'\in\tau_{2}(\ol{\theta},\ell)} 
\!\!\!\Val(\theta';\e,\be_{0}) + 
 \calA_{\ol{\ell}}(\ol{\theta},x_{\ell}) \, 
  \calG_{n+1}(x_{\ell}) \, \partial_{\ell} \!\!\!\sum_{T \in \gotR\gotF_{n}} 
\!\!\! M_{n}(T,x_{\ol{\ell}}) 
  \calG_{n+1}(x_{\ell}) 
  \, \BB_{\ol{\ell}}(\ol{\theta})\\ 
&=\calA_{\ol{\ell}}(\ol{\theta},x_{\ell}) \, 
\calG_{n+1}(x_{\ell})  M_{n}(x_{\ell}) 
\partial_{\be_{0}}\calG_{n}(x_{\ell})\left(1+ 
M_{n}(x_{\ell})\calG_{n+1}(x_{\ell})\right) \BB_{\ol{\ell}}(\ol{\theta}) , 
\end{aligned} 
\end{equation} 
where we have used again that $x_{\ell}=x_{\ol{\ell}}$. Thus, one can reason 
as in (\ref{eq:a.13}), so as to obtain the sum 
of three contributions, as represented in Figure \ref{fig:10}. 
 
\noindent 
\textbf{5.} Finally, consider the case in which $\ell$ is the 
exiting line of a renormalised left-fake cluster, 
$T_{0}$ and the removal of $T_{0}$ and $\ell'_{T_{0}}$ 
creates a self-energy cluster. 
 
Set (for a reason that will become clear later) 
$\theta_{0}=\theta$ and $\ell_{0}=\ell$. 
Then there is a maximal $m\ge1$ such that 
there are $2m$ lines $\ell_{1},\ldots,\ell_{m}$ and 
$\ell'_{1},\ldots\ell'_{m}$, with the following properties: 
 
\noindent 
(i) $\ell_{i}\in \calP(\ell_{\theta_{0}},\ell_{i-1})$, for $i=1,\ldots,m$, \\ 
(ii) $n_{\ell_{i}}=n+i<\max\{p:\Psi_{p}(x_{\ell_{i}})\ne0\}=n+i+1$, for 
$i=0,\ldots,m-1$, while $n_{m}:=n_{\ell_{m}}= 
n+m+\s$, with $\s\in\{0,1\}$,\\ 
(iii) $\nn_{\ell_{i}}\ne\nn_{\ell_{i-1}}$ and the lines preceding 
$\ell_{i}$ but not $\ell_{i-1}$ are on scale $\le n+i-1$, for 
$i=1,\ldots,m$, \\ 
(iv) $\nn_{\ell'_{i}}=\nn_{\ell_{i}}$, for $i=1,\ldots,m$,\\ 
(v) if $m\ge2$, $\ell'_{i}$ is the exiting line of a left-fake cluster 
$T_{i}$, for $i=1,\ldots,m-1$, \\ 
(vi) $\ell'_{i}\prec\ell'_{T_{i-1}}$ and all the lines 
preceding $\ell'_{T_{i-1}}$ but not $\ell'_{i}$ are on scale 
$\le n+i-1$,  for $i=1,\ldots,m$, \\ 
(vii) $n'_{m}:=n_{\ell'_{m}}=n+m+\s'$ with $\s'\in\{0,1\}$.

Note that one cannot have $\s=\s'=1$, otherwise the subgraph between
$\ell_{m}$ and $\ell'_{m}$ would be a self-energy cluster. 
Note also that (ii), (iv) and (v) imply $n_{\ell'_{i}}=n+i$ for $i=1,\ldots,m-1$ 
if $m\ge2$.
Call $S_{i}$ the subgraph between $\ell_{i+1}$ and $\ell_{i}$, and $S'_{i}$ 
the cluster between $\ell'_{T_{i}}$ and $\ell'_{i+1}$ for all 
$i=0,\ldots,m-1$. 
For $i=1,\ldots,m$, call $\theta_{i}$ the renormalised tree obtained from 
$\theta_{0}$ by removing everything between $\ell_{i}$ and the part of  
$\theta_{0}$ preceding $\ell'_{i}$, and note that if $m\ge2$,
properties (i)--(vii) hold for $\theta_{i}$ but with $m-i$ instead of $m$,
for all $i=1,\ldots,m-1$.

\begin{figure}[ht] 
\vskip.3truecm 
\centering 
\ins{018pt}{-23pt}{$\theta_{0}=$} 
\ins{072pt}{-16pt}{$n\!+\!1$} 
\ins{080pt}{-32pt}{$\ell_{1}$} 
\ins{125pt}{-24pt}{$\le n$} 
\ins{150pt}{-46pt}{$S_{0}$} 
\ins{177pt}{-17pt}{$n$} 
\ins{177pt}{-32pt}{$\ell_{0}$} 
\ins{226pt}{-25pt}{$n$} 
\ins{250pt}{-46pt}{$T_{0}$} 
\ins{268pt}{-16pt}{$n\!+\!1$} 
\ins{274pt}{-32pt}{$\ell_{T_{0}}'$} 
\ins{314pt}{-24pt}{$\le n$} 
\ins{344pt}{-46pt}{$S'_{0}$} 
\ins{370pt}{-32pt}{$\ell'_{1}$} 
\ins{364pt}{-16pt}{$n\!+\!1$} 
\ins{018pt}{-82pt}{$\theta_{1}=$} 
\ins{072pt}{-75pt}{$n\!+\!1$} 
\ins{080pt}{-92pt}{$\ell_{1}$} 
\ins{226pt}{-88pt}{$R_{1}$} 
\ins{226pt}{-129pt}{$n$} 
\ins{226pt}{-143pt}{$\ell_{0}$} 
\ins{177pt}{-133pt}{$\le n$} 
\ins{274pt}{-133pt}{$\le n$} 
\ins{253pt}{-155pt}{$S'_{0}$} 
\ins{203pt}{-155pt}{$S_{0}$} 
\includegraphics[width=5.2in]{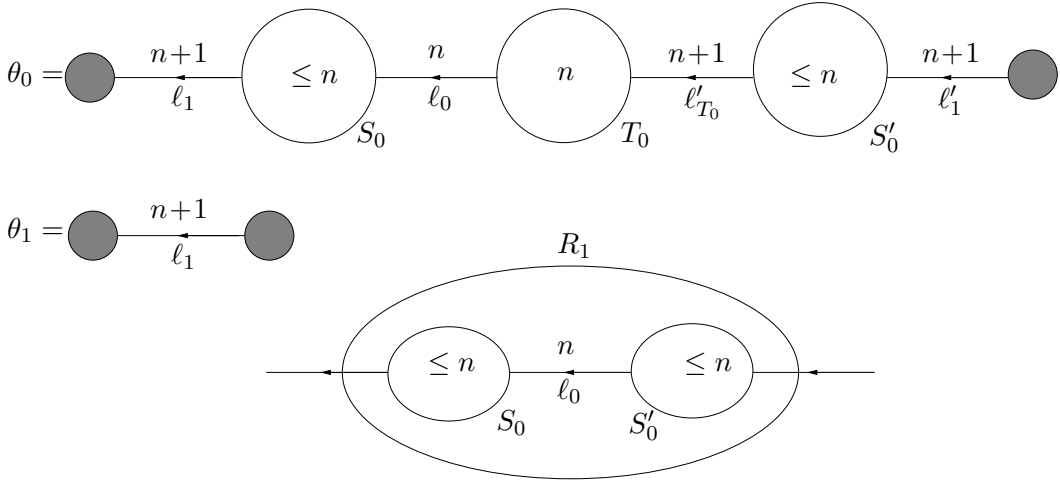} 
\caption{The renormalised trees $\theta_{0}$ and $\theta_{1}$ and the 
self-energy cluster $R_{1}$ in case 5 
with $m=1$ and $\s=\s'=0$. Note that the set $S_{0}'$ is a 
cluster, but not a self-energy cluster.} 
\label{fig:11} 
\end{figure}

For $i=1,\ldots,m$, call $R_{i}$ the self-energy 
cluster obtained from the subgraph of $\theta_{i-1}$ between $\ell_{i}$ and
$\ell'_{i}$, by removing the left-fake cluster $T_{i-1}$ together with
$\ell'_{T_{i}}$. Note that $L(R_{i})=L(S_{i-1})\cup\{\ell_{i-1}\}\cup 
L(S'_{i-1})$ and $N(R_{i})=N(S_{i-1})\cup N(S'_{i-1})$; 
see Figure  \ref{fig:11}.
 
For $i=0,\ldots,m-1$, given $\ell',\ell\in L(\theta_{i})$, with 
$\ell'\prec\ell$, call $\calP^{(i)}(\ell,\ell')$ the path of lines in 
$\theta_{i}$ connecting 
$\ell'$ to $\ell$ (hence $\calP^{(i)}(\ell,\ell')=\calP(\ell,\ell') 
\cap L(\theta_{i})$). For any $i=0,\ldots,m-1$ and any $\ell\in\calP^{(i)} 
(\ell_{i},\ell'_{m})$, let 
$\tau_{3}(\theta_{i},\ell)$ be the set of all renormalised trees 
which can be obtained from $\theta_{i}$ by replacing each left-fake cluster 
preceding $\ell$ but not $\ell'_{m}$ with all possible left-fake 
clusters. Set also $\tau_{3}(\theta_{m-1},\ell_{m}')=\theta_{m-1}$. 
 
Note that 
\begin{equation}\label{eq:a.18} 
\begin{aligned} 
&\calA_{\ell_{m}}(\theta_{m},x_{\ell_{m}}) 
  \calG_{n_{m}}(x_{\ell_{m}}) \Val(S_{m-1})= 
\calA_{\ell_{m-1}}(\theta_{m-1},x_{\ell_{m-1}}),\\ 
&\Val(S'_{m-1}) \calG_{n'_{m}}(x_{\ell_{m}}) 
\BB_{\ell_{m}}(\theta_{m}) 
=\BB_{\ell'_{T_{m-1}}}(\theta_{m-1}), 
\end{aligned} 
\end{equation} 
and one among cases 1--4 holds for $\ell_{m}\in L(\theta_{m})$ 
so that we can consider the contribution to $\partial_{\ell_{m}}\Val 
(\theta_{m};\e,\be_{0})$ (together with other contributions as in 3 and 4 
if necessary) given by -- see (\ref{eq:a.10}), (\ref{eq:a.11}) and
(\ref{eq:a.13}) -- 
$$ 
\calA_{\ell_{m}}(\theta_{m},x_{\ell_{m}})\calG_{n_{m}}(x_{\ell_{m}}) 
\partial_{\ell_{m-1}}\Val_{R_{m}}(x_{\ell_{m}})\calG_{n'_{m}}(x_{\ell_{m}}) 
\BB_{\ell_{m}}(\theta_{m}). 
$$ 
Then one has 
\begin{equation}\label{eq:a.19} 
\begin{aligned} 
&\calA_{\ell_{m}}(\theta_{m},x_{\ell_{m}})\calG_{n_{m}}(x_{\ell_{m}}) 
\partial_{\ell_{m-1}}\Val_{R_{m}}(x_{\ell_{m}})\calG_{n'_{m}}(x_{\ell_{m}}) 
\BB_{\ell_{m}}(\theta_{m})+ 
\partial_{\ell_{m-1}}\!\!\!\!\!\!\!\!\!\!\! 
\sum_{\theta'\in \tau_{3}(\theta_{m-1},\ell_{m-1})} 
\!\!\!\!\!\!\!\!\!\!\! 
\Val(\theta';\e,\be_{0})\\ 
&=\calA_{\ell_{m-1}}(\theta_{m-1},x_{\ell_{m-1}})\partial_{\be_{0}} 
\calG_{n+m-1}(x_{\ell_{m-1}}) 
\left( 1+M_{n+m-1}(x_{\ell_{m-1}})\calG_{n+m}(x_{\ell_{m-1}}) \right) 
\\ &\qquad\times 
\BB_{\ell'_{T_{m-1}}} (\theta_{m-1}), 
\end{aligned} 
\end{equation} 
and hence we obtain, reasoning as in (\ref{eq:a.13}), 
\begin{equation}\label{eq:a.20} 
\begin{aligned} 
&\calA_{\ell_{m-1}}(\theta_{m-1},x_{\ell_{m-1}}) 
 \calG_{n+m-1}(x_{\ell_{m-1}})\partial_{\be_{0}}\MM_{n+m-2}(x_{\ell_{m-1}}) 
\calG_{n+m-1}(x_{\ell_{m-1}}) 
\BB_{\ell'_{T_{m-1}}}(\theta_{m-1})\\ 
&\qquad+\calA_{\ell_{m-1}}(\theta_{m-1},x_{\ell_{m-1}}) 
 \calG_{n+m-1}(x_{\ell_{m-1}})\partial_{\be_{0}}\MM_{n+m-2}(x_{\ell_{m-1}}) 
\calG_{n+m}(x_{\ell_{m-1}})\\ 
&\qquad\qquad\times 
\BB_{\ell'_{T_{m-1}}}(\theta_{m-1})\\ 
&\qquad+\calA_{\ell_{m-1}}(\theta_{m-1},x_{\ell_{m-1}}) 
 \calG_{n+m-1}(x_{\ell_{m-1}})\partial_{\be_{0}}\MM_{n+m-2}(x_{\ell_{m-1}}) 
\calG_{n+m-1}(x_{\ell_{m-1}})\\ 
&\qquad\qquad\times 
M_{n+m-1}(x_{\ell_{m-1}}) \calG_{n+m}(x_{\ell_{m-1}}) 
\BB_{\ell'_{T_{m-1}}}(\theta_{m-1}). 
\end{aligned} 
\end{equation} 

Then, for $i=m-1,\ldots,1$ we recursively reason as follows. Set 
$$ 
\BB_{\ell'_{T_{i}}}(\tau_{3}(\theta_{i},\ell'_{i+1})):= 
\sum_{\theta'\in\tau_{3}(\theta_{i},\ell'_{i+1})} 
\BB_{\ell'_{T_{i}}}(\theta') 
$$ 
and note that 
\begin{equation}\label{eq:a.21} 
\begin{aligned} 
&\calA_{\ell_{i}}(\theta_{i},x_{\ell_{i}}) 
  \calG_{n+i}(x_{\ell_{i}}) \Val(S_{i-1})= 
\calA_{\ell_{i-1}}(\theta_{i-1},x_{\ell_{i-1}}),\\ 
&\Val(S'_{i-1}) \calG_{n+i}(x_{\ell_{i}}) 
M_{n+i}(x_{\ell_{i}}) \calG_{n+i+1}(x_{\ell_{i}}) 
\BB_{\ell'_{T_{i}}}(\tau_{3}(\theta_{i},\ell'_{i+1}) )
=\BB_{\ell'_{T_{i-1}}}(\tau_{3}(\theta_{i-1},\ell'_{i})). 
\end{aligned} 
\end{equation} 
Consider the contribution 
\begin{equation}\label{eq:a.22} 
\calA_{\ell_{i}}(\theta_{i},x_{\ell_{i}}) 
  \calG_{n+i}(x_{\ell_{i}}) 
  \partial_{\ell_{i-1}}\Val_{R_{i}}(x_{\ell_{i}}) 
  \calG_{n+i}(x_{\ell_{i}}) 
M_{n+i}(x_{\ell_{i}})  \calG_{n+i+1}(x_{\ell_{i}}) 
\BB_{\ell'_{T_{i}}}(\tau_{3}(\theta_{i},\ell'_{i+1}) )
\end{equation} 
obtained at the $(i+1)$-th step of the recursion. 
By (\ref{eq:a.21}) one has (see Figure \ref{fig:12}) 
\begin{equation}\label{eq:a.23} 
\begin{aligned} 
&\calA_{\ell_{i}}(\theta_{i},x_{\ell_{i}}) 
  \calG_{n+i}(x_{\ell_{i}}) 
  \partial_{\ell_{i-1}}\Val_{R_{i}}(x_{\ell_{i}}) 
  \calG_{n+i}(x_{\ell_{i}}) 
M_{n+i}(x_{\ell_{i}})  \calG_{n+i+1}(x_{\ell_{i}}) 
\BB_{\ell'_{T_{i}}}(\tau_{3}(\theta_{i},\ell'_{i+1}))\\ 
&\qquad+ 
\partial_{\ell_{i-1}}\!\!\!\!\!\!\!\!\!\!\! 
\sum_{\theta'\in \tau_{3}(\theta_{i-1},\ell_{i-1})} 
\!\!\!\!\!\!\!\!\!\!\! 
\Val(\theta';\e,\be_{0}) 
=\calA_{\ell_{i-1}}(\theta_{i-1},x_{\ell_{i-1}}) \, 
\partial_{\be_{0}}\calG_{n+i-1}(x_{\ell_{i-1}}) 
\\ & \times \left(1+ 
M_{n+i-1}(x_{\ell_{i-1}})\calG_{n+i}(x_{\ell_{i-1}})\right) 
 \BB_{\ell'_{T_{i-1}}}(\tau_{3}(\theta_{i-1},\ell'_{i})), 
\end{aligned} 
\end{equation} 
which produces, as in (\ref{eq:a.20}), the contribution 
\begin{equation}\label{eq:a.24} 
\begin{aligned} 
\calA_{\ell_{i-1}}(\theta_{i-1},x_{\ell_{i-1}}) 
  & \calG_{n+i-1}(x_{\ell_{i-1}}) 
  \partial_{\ell_{i-2}}\Val_{R_{i-1}}(x_{\ell_{i-1}}) 
  \calG_{n+i-1}(x_{\ell_{i-1}})
\\ &\times 
  M_{n+i-1}(x_{\ell_{i-1}})\calG_{n+i}(x_{\ell_{i-1}})
  \BB_{\ell'_{T_{i-1}}}(\tau_{3}(\theta_{i-1},\ell'_{i})). 
\end{aligned} 
\end{equation} 
%
 
Hence we can proceed recursively from $\theta_{m}$ up to $\theta_{0}$, 
until we obtain 
\begin{equation}\label{eq:a.25} 
\begin{aligned} 
&\calA_{\ell_{0}}(\theta_{0},x_{\ell_{0}}) 
 \calG_{n}(x_{\ell_{0}})\partial_{\be_{0}}\MM_{n-1}(x_{\ell_{0}}) 
 \calG_{n}(x_{\ell_{0}})  
  \BB_{\ell'_{T_{0}}}(\tau_{3}(\theta_{0},\ell'_{1}))\\ 
&\qquad+\calA_{\ell_{0}}(\theta_{0},x_{\ell_{0}}) 
 \calG_{n}(x_{\ell_{0}})\partial_{\be_{0}}\MM_{n-1}(x_{\ell-{0}}) 
 \calG_{n+1}(x_{\ell_{0}})  
 \BB_{\ell'_{T_{0}}}(\tau_{3}(\theta_{0},\ell'_{1}))\\ 
&\qquad+\calA_{\ell_{0}}(\theta_{0},x_{\ell_{0}}) 
 \calG_{n}(x_{\ell_{0}})\partial_{\be_{0}}\MM_{n-1}(x_{\ell_{0}}) 
 \calG_{n}(x_{\ell_{0}})M_{n}(x_{\ell_{0}}) \calG_{n+1}(x_{\ell_{0}}) 
 \, \BB_{\ell'_{T_{0}}}(\tau_{3}(\theta_{0},\ell'_{1})). 
\end{aligned} 
\end{equation} 
Once again, if we split $\partial_{\be_{0}}=\partial_{N}+\partial_{L}$, 
all the terms with $\partial_{N}\MM_{n-1}$ are contributions to 
$\MM_{\io}(0)$. 
 
\noindent 
\textbf{6.} We are left with the derivatives $\partial_{L}M_{q}(x;\e,\be_{0})$, $q\le n$, 
when the differentiated propagator is not one of those used along 
the cases 4 or 5; see for instance (\ref{eq:a.16}), (\ref{eq:a.19}) and 
(\ref{eq:a.23}). One can reason 
as in the case $\partial_{L}\Val(\theta;\e,\be_{0})$, by studying the 
derivatives 
$\partial_{\ell}\Val_{T}(x_{\ell};\e,\be_{0})$ and proceed iteratively 
along the lines of cases 1 to 5 above, until 
only lines on scales $0$ are left. In that case the derivatives 
$\partial_{\be_{0}}\calG_{0}(x_{\ell};\e,\be_{0})$ 
produce derivatives $\partial_{\be_{0}}M_{-1}(x;\e,\be_{0})=\e 
\partial_{\be_{0}}^2F_{\vzero}(\beta_{0})$ 
(see Remarks \ref{rmk:3.5bis} and \ref{rmk:3.1ter}). Therefore, for $n=-1$, 
in the splitting (\ref{eq:a.9}), there are no terms with the 
derivatives $\partial_{\ell}$, and the derivatives $\partial_{v}$ can 
be interpreted as said in Remark \ref{rmk:a.1}. It is also easy 
to realize that, by construction, each contribution to $\MM_{\io}(0;\e,\be_{0})$ 
appears as one term among those considered in the discussion above. 
Hence the assertion follows. 
\EP 

\begin{rmk}\label{rmk:a.5}
\emph{
If we used a sharp scale decomposition instead of the $C^{\io}$ one,
the proof above would be much more easier.
More precisely, if we defined the (discontinuous) function
\begin{equation}\nonumber
\chi(x):=\left\{\begin{aligned} 
&1,\qquad |x| \le 1, \\ 
&0,\qquad |x| > 1, 
\end{aligned}\right. 
\end{equation} 
and consequently changed the definitions of $\psi$, and $\chi_{n},\psi_{n}$
and $\Psi_{n}$ for $n\ge0$, we could reduce the proof of Lemma \ref{lem:4.8}
to (iterations of) case 1.
Moreover in such a case, setting
$$
G^{\RR}_{n}(\e,\be_{0})=\sum_{k\ge0}\e^{k}G^{[k]}_{n}(\e,\be_{0}),
\qquad
G^{[k]}_{n}(\e,\be_{0})=
\sum_{\theta\in\Theta^{\RR}_{k+1,\vzero,n}}\Val(\theta,\e,\be_{0}),
$$
with $\Theta^{\RR}_{k,\nn,n}=\{\theta\in\Theta^{\RR}_{k,\nn,n}\,:\,n_{\ell}\le
n \mbox{ for all }\ell\in L(\theta)\}$, we would obtain the stronger
identity
$$
\MM_{n}(0;\e,\be_{0})=\e\partial_{\be_{0}}G^{\RR}_{n}(\e,\be_{0}),
$$
for all $n\ge-1$.
On the other hand, the bound (\ref{eq:4.3b}) in Lemma \ref{lem:4.6}
would be no longer true because of the derivative $\partial_{x}\Psi_{n}$,
so that further work would be however needed;
see for instance \cite{GBG} where a sharp scale decomposition is used
for the standard KAM theorem and $\oo$ satisfying the standard Diophantine
condition.
}
\end{rmk}
 

\end{document}